\documentclass[fleqn]{mat01}
\usepackage{times,mathtimy,amssymb,epsfig}%krepsf,psfig}
\begin{document}

\setcounter{page}{123}
\firstpage{123}

\newcommand{\myfavouritearrow}{$\leftrightsquigarrow$}

\newtheorem{theore}{Theorem}
\renewcommand\thetheore{\arabic{section}.\arabic{theore}}
\newtheorem{theor}[theore]{\bf Theorem}
\def\notation{\trivlist\item[\hskip\labelsep{{\it Notation.}}]}
\def\teo{\trivlist\item[\hskip\labelsep{\bf Theorem}]}
\def\potc{\trivlist\item[\hskip\labelsep{{\it Proof of the Claim.}}]}
\newtheorem{theorr}{Theorem}
\renewcommand\thetheorr{\arabic{theorr}}
\newtheorem{therr}[theorr]{\bf Theorem}
\newtheorem{lem}{Lemma}
\newtheorem{propo}{\rm PROPOSITION}
\newtheorem{rem}{Remark}
\newtheorem{caa}{Case}
\renewcommand\thecaa{(\alph{caa})}
\newtheorem{case}[caa]{Case}
\newtheorem{coro}{\rm COROLLARY}
\newcommand{\A}{{\mathbb A}}
\newcommand{\C}{{\mathbb C}}
\newcommand{\V}{{\mathbb V}}
\newcommand{\Z}{{\mathbb Z}}
\def\o{{\cal O}}
\def\L{{\cal L}}
\def\P{{\Bbb P}}
\def\R{{\Bbb R}}
\def\G{{\cal G}}
\def\E{{\cal E}}
\def\D{{\cal D}}

\font\zz=msam10 at 10pt
\def\BBox{\mbox{\zz{\char'245}}}
\def\Box{\mbox{\zz{\char'244}}}

\makeatletter
\newcommand\abb{\@setfontsize\abb{7.5}{9}}
\def\artpath#1{\def\@artpath{#1}}
\makeatother
\artpath{c:/prema/pm}

\title{Twisted holomorphic forms on generalized flag varieties}

\markboth{K Paramasamy}{Twisted holomorphic forms on generalized flag varieties}

\author{K PARAMASAMY}

\address{Chennai Mathematical Institute, Chennai~600~017, India\\
\noindent E-mail: paramas@cmi.ac.in}

\volume{114}

\mon{May}

\parts{2}

\Date{MS received 22 September 2003; revised 13 April 2004}

\begin{abstract}
In this paper we prove some vanishing theorems for the twisted Dolbeault
cohomology of the complete flag varieties associated to a simple, simply
connected algebraic group.
\end{abstract}


\maketitle

\section{Introduction}

Let $G$ be a simple, simply connected algebraic group defined over an
algebraically closed field $k$ of characteristic zero. Fix a maximal
torus $T$. Let $W = \hbox{Nor}_G(T)/\hbox{Cent}_G(T)$ be the Weyl group
and $X(T)=\hbox{Hom}(T, GL_1(k))$ be the group of characters. The main
results of this article are as follows:

\begin{teo} {\bf (Theorem 1).}\ \ {\it Let $G$ be a simple simply
connected algebraic group over an algebraically closed field $k$ of
characteristic zero. Let $p \geq 0.$ Let $\lambda$ be a dominant weight
such that for each $\mu \in \Phi^-_p$ either $\lambda + \mu$ is dominant
or $\lambda + \mu + \rho$ is singular. Then $H^{p,q}(G/B,
\L ({\lambda})) = 0$ for all $q\geq 1.$}
\end{teo}\vspace{.6pc}

Let $\lambda= \sum _i n_i \omega _i$.  For given $p \geq 0$, in
Proposition~2 we give conditions on $n_i$'s so that $\lambda$
satisfies the hypothesis of Theorem~1.

For each simple root $\alpha$, we define the {\em Coxeter number of
$\alpha$}, denoted by $ h_{\alpha}$, as the number $\sum_{\gamma \in
\Phi_{\alpha}^+}(\gamma, \alpha^{\vee})$, where
$\Phi^+_{\alpha}=\{\gamma \in \Phi^+ | (\gamma, \alpha^{\vee})>0\}$. We
prove in Theorem~9 that $h_{\alpha} \leq h$, and for a shorter simple
root $\beta, h_{\beta}=h$, where $h$ is the Coxeter number\break of $G$.

As a corollary of Theorem~1 and Proposition~2, we get the following
vanishing theorem:
\begin{equation*}
H^{p,q}(G/B, \L_{\lambda})=0,\ \  q \geq 1,\ \  p \geq 0,\ \ \hbox{if}\ \ n_k \geq
h_{\alpha_k}-1,\ \ \forall k.
\end{equation*}

Let rank$(G) \geq 2$ and $\dim(G/B)=d.$ In Theorem~12 we give strictly
dominant weight $\lambda$ such that $H^{d-1,1}(G/B,
\L({\lambda}))\neq 0.$

We also derive the following corollaries:

\begin{enumerate}
\renewcommand{\labelenumi}{\arabic{enumi}.}
\item We show that the Bott vanishing property fails to hold for
generalised flag varieties.
\item $G/B$ is not a toric variety.
\item $G/B$ cannot be degenerated to a smooth toric variety in such a
manner that ample cone degenerates to ample cone.\vspace{-.5pc}
\end{enumerate}

We recall some standard facts about algebraic groups. A standard
reference for this material is \cite{jant}. For a finite dimensional
$T$-module $V$ and ${\chi \in X(T)}$, set $V_{\chi}= \{v \in V | tv
= \chi (t)v, \ \forall t \in T\}.$ The finitely many $\chi \in X(T)$
such that $V_{\chi} \not = 0$ are called weights of the $T$-module $V$
and $V_{\chi}$ the weight spaces. Then $V$ decomposes into a direct sum
of weight spaces, i.e. $V = \bigoplus_{\chi_i} V_{\chi _i}$, where
${\chi _i}$ are weights.

Let ${\frak g}$ and ${\frak h}$ be the Lie algebras of $G$ and $T$
respectively. Let $T$ act on ${\frak g}$ via the adjoint representation.
The non-zero weights of this representation are called roots and we
denote the set of all roots by ${\Phi}$. The weight space corresponding
to the zero weight is ${\frak h}$. We identify $X(T)$ canonically as a
subset of $X(T)\otimes{\R}$. Let $\Delta$ denote the set
of all simple roots, that is a subset of ${\Phi}$ such that (i)
$\Delta$ is a basis of $X(T) \otimes {\R}$, (ii) for each
root ${\alpha \in \Phi}$ there exist integers ${n_i}$ of like sign such
that ${\alpha = \sum n_i \alpha_i}, \ \alpha_i \in \Delta$. The
number of elements in $\Delta$ is called the rank of the group
$G$. The set of all roots ${\alpha \in \Phi}$ for which the $n_i$ are
not negative is denoted by ${\Phi^+}$ and an element of ${\Phi^+}$ is
called a positive root. ${\Phi^-: = \Phi -\Phi^+}$ is the set of
negative roots. Let $ \frak{n} =\bigoplus _{\alpha \in \Phi^+}{\frak g}
_{\alpha}$ and $ \frak{n}^- =\bigoplus _{\alpha \in \Phi^-}{\frak g}
_{\alpha}$. Let $B$ be the Borel subgroup of $G$ such that the Lie
algebra ${\frak b}$ of $B$\break is ${\frak h} \oplus \frak{n^-}$.

There is a natural faithful action of $W$ on $X(T)\otimes {\R}$ and
there exists a $W$-invariant bilinear form $(,)$ on $X(T) \otimes
{\R}$. For $\alpha \in \Delta$, we denote $2\alpha/(\alpha, \alpha)$ by
$\alpha^{\vee}.$ For each $\alpha \in \Phi$, we have reflection
$s_\alpha \in W$ defined by $s_\alpha(\psi) = \psi-(\psi,
\alpha^{\vee})\alpha$, $\forall \psi \in X(T)\otimes{\R}$. The set
$\{s_\alpha, \alpha \in \Delta\}$ generates $W$. For each
$\alpha_i \in \Delta$ there exists a weight $\omega_i \in X(T)$
such that $(\omega_i, \alpha_j^{\vee}) = \delta_{ij}$ (Here $\delta$ is
the Kr\"{o}necker delta.) These $\omega _i$ are called fundamental weights.
The action of $W$ permutes the roots. For $w \in W$ we define the length
of $w$, denoted by $l(w),$ to be the number of positive roots moved by
$w$ to negative\break roots.

Following Jantzen \cite{jant}, we denote the sheaf associated to a
$B$-module $V$ by $\L(V) = \L_{G/B}(V)$ (i.e. the sheaf of sections of
the associated vector bundle on $G/B$). When $V=k_{\lambda}$, we denote
$\L(k_{\lambda})$ by $\L(\lambda)$. The tangent sheaf on $G/B$
(\cite{jant},\ p.~229) is $\L(\frak{g}/\frak{b} )$. The $T$ weights on
$\frak{g}/\frak{b}$ are $\Phi^+$.

We say a weight $\lambda = \sum n_k \omega_k$ is dominant (strictly
dominant) if $n_k \geq 0 ~ (n_k>0).$ We denote the set of dominant
weight by $X(T)^+$. For $w \in W, ~ \lambda \in X(T)$ we have the `dot'
action given by $w \cdot \lambda := w(\lambda+\rho)-\rho,$ where $\rho$
is half sum of positive roots. Now we recall the Borel--Weil--Bott
theorem. Let $\lambda \in X(T)$. If $\lambda \in X(T)^+$ then $H^i(G/B,
\L({\lambda})) \neq 0$ only when $i=0$, and $\{H^0(G/B,
\L({\lambda}))~ \lambda \in X(T)^+\}$ is the set of all simple
$G$-modules. If there is a $w \in W$ such that $w \cdot \lambda \in
X(T)^+$ (such a $w$ will be unique) then $H^i(G/B, \L(\lambda))
\neq 0$ only when $i=l(w)$ and $H^{l(w)} (G/B,
\L(\lambda))\simeq H^0(G/B, \L({w\cdot \lambda}))$
and moreover if there is no such $w$ (equivalently if $\lambda + \rho$
is singular, i.e. $(\lambda + \rho, \alpha ^{\vee}) =0$ for some root
$\alpha \in \Phi$) then $H^i(G/B, \L(\lambda))=0, ~ \forall i$.

Recall that in \cite{bott}, Bott proved the following vanishing theorem
for the complex projective space: $H^*({\Bbb P}^n, \Omega_{{\Bbb
P}^n}^q\otimes \o (k))=0, ~ 1\leq k \leq q.$ We say that a
smooth projective variety $X$ satisfies the Bott vanishing
property if the twisted holomorphic forms or Dolbeault cohomology
$H^{p,q}(X,\L)=H^q(X,\Omega^p(\L))$ $=0$ for all
ample line bundle $L,~ \forall q>0~{\rm and}~ \forall p\geq 0,$ where
$\Omega ^p$ stands for sheaf associated to the vector bundle $\Lambda ^p
T ^*_X$ and $\Omega^p(\L) = \Omega^p \otimes \L.$

We use the partial order `$\geq$' in $X(T)$, for $\lambda_1,~ \lambda_2
\in X(T)$, we say $\lambda_1 \geq \lambda_2$ if $\lambda_1 -\lambda_2 =$
sum of some simple roots. For a simple algebraic group, the Coxeter
number is defined as $\frac {\dim(G)}{{\rm rank}(G)}-1$. For the
ordering of the simple roots we follow the following diagrams
\cite{serre}. In the following diagram the arrows are pointing towards
shorter simple root. If there is no arrow it means that the adjacent roots
are of the same length.

\begin{figure}[h]
%\begin{center}
%\parbox{0.55\textwidth}{\epsfxsize=0.55\textwidth
\centerline{\epsfxsize=9cm\epsfbox{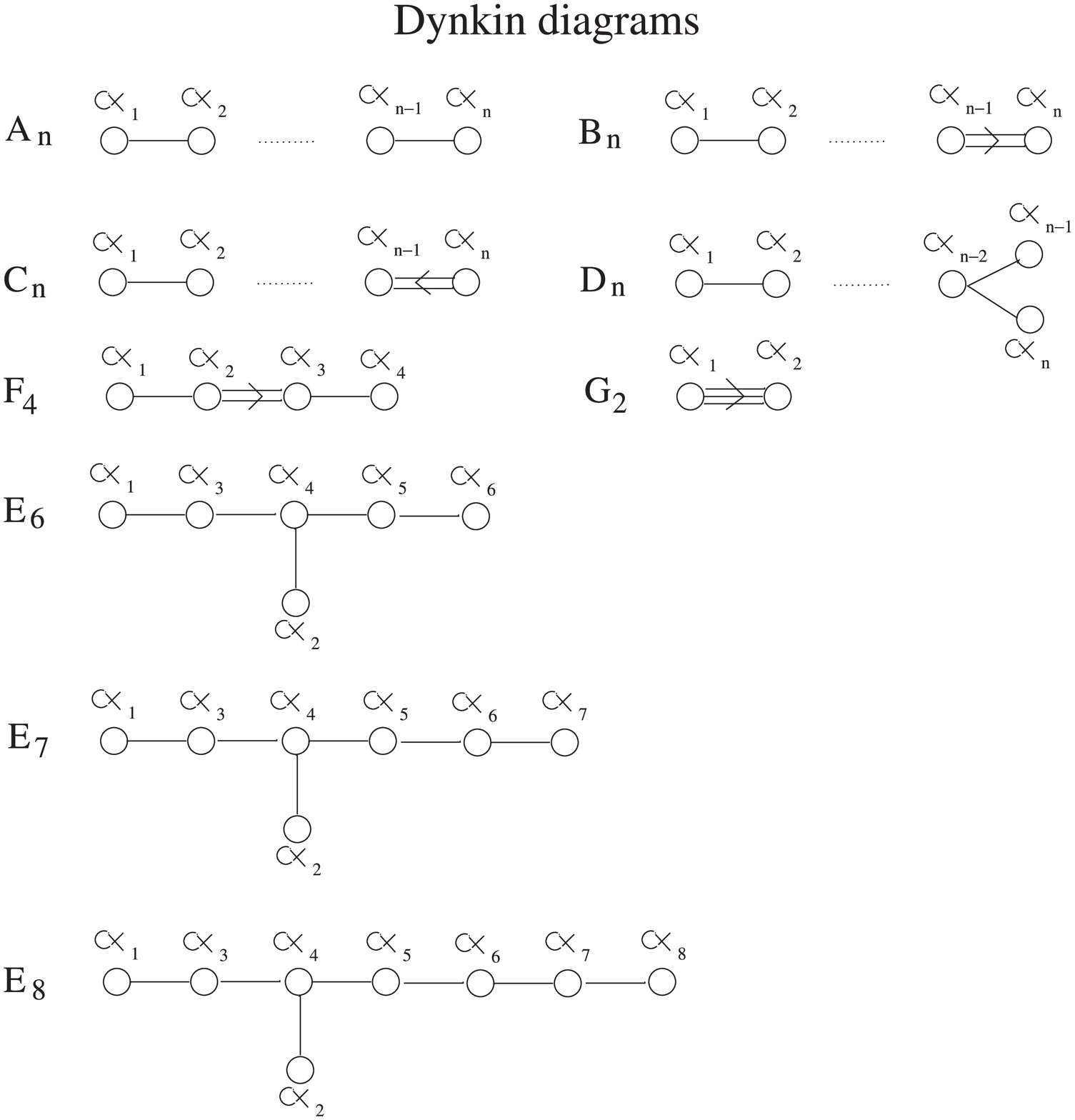}}\vspace{.8pc}
%\end{center}
\end{figure}

\subsection{\it New notations}

Now we introduce some new notations. For any subset $S$ of $X(T)$ and
$j\geq 1$, $S_j := \{ \nu \in X(T) | \nu= \sum_{i=1}^j \lambda_i,
\lambda_i \in S$ and $\lambda_i$'s are $\hbox{distinct}\}$. For any finite set
$S$ of $X(T)$, we define an $|S| \times {\rm rank}(G) $-matrix
associated to $S$ as follows: We represent $\lambda = \sum n_k \omega_k
\in S$ as tuples $(n_1, n_2, \ldots )$, and write the elements of $S$
one below the other to get a $|S| \times {\rm rank}(G)$ matrix. We call
the matrix corresponding to $\Phi ^+$ the positive roots matrix. For
$\alpha \in \Delta, \Phi^+_{\alpha}: = \{ \beta \in \Phi^+|
(\beta, \alpha^{\vee})>0\}.$ We define the Coxeter number of
$\alpha = h_{\alpha}= \sum_{\beta \in \Phi^+_{\alpha}}(\beta,
\alpha^{\vee}).$

\section{Vanishing theorems}

In this section, we study conditions on dominant weights $\lambda$ so
that $H^{p,q}(G/B,\L({\lambda})) =0$ for all $q \geq 1$ and
$p \geq 0$.

\begin{therr}[\!]
Let $p \geq 0.$ Let $\lambda$ be a dominant weight such that for each
$\mu \in \Phi^-_p $ either $\lambda + \mu$ is dominant or $\lambda + \mu
+ \rho$ is singular. Then $H^{p,q}(G/B, \L({\lambda})) = 0$ for all
$q\geq 1.$
 \end{therr}

\begin{proof}
The proof follows from the following observations:
\begin{enumerate}
\renewcommand{\labelenumi}{\arabic{enumi}.}
\leftskip -.2pc
\item By the hypothesis and by the Borel--Weil--Bott theorem, for any
weight $\nu$ of $\Lambda ^p \frak{n^-} \otimes k_{\lambda}$, we have
$H^q(G/B, \L({\nu}))=0$ for all $q\geq1$.

\item Any lowest weight $\lambda_1$ of $\Lambda ^p \frak{n^-} \otimes
k_{\lambda }$ will give a one-dimensional $B$-submodule $k_{\lambda_1}$.
Let $V_1$ be $V/k_{\lambda _1}$. Then we will have a short exact
sequence of $B$-modules; $0\longrightarrow k_{\lambda_1}
\longrightarrow V \longrightarrow V_1 \longrightarrow 0.$

\item The above short exact sequence induces the following long exact
sequence:
\begin{align*}
\hskip -.55cm 0 &\rightarrow H^0(G/B,\L({\lambda_1})) \rightarrow H^{p,0}(G/B,
  \L({\lambda}))\rightarrow H^0(G/B, \L(V_1))\\[.2pc]
\hskip -.55cm &\rightarrow H^1(G/B, \L({\lambda_1}))
\rightarrow H^{p,1}(G/B, \L({\lambda}))
\rightarrow H^1(G/B,\L(V_1))\\[.2pc]
\hskip -.55cm &\rightarrow H^2(G/B, \L({\lambda_1}))
\rightarrow H^{p,2}(G/B, \L({\lambda})) \rightarrow
   \cdots \cdots \cdots
\end{align*}
\item By $1$ and $3$, we have $H^{p,q}(G/B,
\L_{\lambda})\simeq H^q(G/B, \L(V_1))$ for all $q\geq 1$. Now we repeat
steps $2$ and $3$ with the $B$-module $V_1$. Let $\lambda_2$ be the
lowest weight of $V_1$ and $V_2=V_1/k_{\lambda_1}$. As before we get
$H^q(G/B, \L(V_1) \simeq H^q(G/B, \L(V_2)$ for all $q\geq1$. After a
finite number of steps we end up with a two-dimensional $B$-module
$V_j$ such that $H^q(G/B, \L(V_j))=0$ for all $q\geq 1$. This proves the
theorem.\hfill $\BBox\ \,$
\end{enumerate}\vspace{-2pc}
\end{proof}

We now make explicit computations to give more delicate bound on the
$\lambda$ for vanishing results. For this we use the classification of
simple algebraic groups.

\setcounter{propo}{1}
\begin{propo}$\left.\right.$
\begin{enumerate}
\renewcommand\labelenumi{\rm (\Alph{enumi})}
\leftskip .25pc
\item Let $G$ be a simply connected algebraic group of type $A_n$. Let
$\lambda = \sum n_k \omega_k.$ Then
\begin{equation*}
\hskip -.55cm H^{p,q}(G/B, \L({\lambda})) =0, \forall ~q\geq1,
\end{equation*}
if
\begin{equation*}
\hskip -.55cm \begin{cases}
n_k \geq p &\hbox{when} \hskip .3cm 0\leq p \leq n-1\\
 &\\[-.7pc]
n_k \geq n &\hbox{when} \hskip .3cm n \leq p \leq
\displaystyle\frac{n^2-n+2}{2}\\
 &\\[-.7pc]
n_k \geq \displaystyle\frac{n^2+n+2-2p}{2} &\hbox{when} \hskip
.3cm\displaystyle\frac{n^2-n+4}{2} \leq p \leq \frac{n(n+1)}{2}
\end{cases}.
\end{equation*}
\item Let $G$ be a simply connected algebraic group of type $B_n$. Let
$\lambda=\sum n_k \omega_k.$ Then
\begin{equation*}
\hskip -.55cm H^{p,q}(G/B, \L({\lambda})) =0, \forall ~q\geq1,
\end{equation*}
if
\begin{equation*}
\hskip -.55cm \begin{cases}
n_k \geq p, k \not = n; n_n \geq 2p-1 &\hbox{when} \hskip .3cm 0\leq
p\leq n-1 \\
 &\\[-.7pc]
n_k \geq p, k \not = n;n_n \geq 2n-1 &\hbox{when} \hskip .3cm n\leq
p\leq 2n-3 \\
 &\\[-.7pc]
n_k \geq 2n-2, k \not = n;n_n \geq 2n-1 &\hbox{when} \hskip .3cm
2n-2\leq p \leq n^2-2n+3\\
 &\\[-.7pc]
n_k \geq n^2+1-p,   &\hbox{when} \hskip .3cm n^2-2n+4 \leq p\\
 &\\[-.7pc]
~~k \not= n;n_n \geq 2n-1 &\hskip 1cm \leq n^2-n+1\\
 &\\[-.7pc]
n_k \geq n^2+1-p, &\hbox{when} \hskip .3cm n^2-n+2 \leq p \leq n^2\\
 &\\[-.7pc]
~~  k \not = n;n_n \geq 2n^2+1-2p
\end{cases}.
\end{equation*}
\item Let $G$ be a simply connected algebraic group of type $C_n$. Let
$\lambda=\sum n_k \omega_k$. Then
\begin{equation*}
\hskip -.55cm H^{p,q}(G/B, \L({\lambda})) =0, \forall ~q\geq1,
\end{equation*}
if
\begin{equation*}
\hskip -.55cm \begin{cases}
n_k \geq 2p-1, k \not = n; n_n \geq p &\hbox{when} \hskip .3cm p =1,2\\[.2pc]
n_k \geq p+1,  k \not = n; n_n \geq p &\hbox{when} \hskip .3cm 3\leq p\leq n-1 \\[.2pc]
n_k \geq p+1,  k \not = n; n_n \geq n &\hbox{when} \hskip .3cm n\leq p\leq 2n-3 \\[.2pc]
n_k \geq 2n-1, k \not = n; n_n \geq n &\hbox{when} \hskip .3cm 2n-2\leq p \leq n^2-2n+3\\[.2pc]
n_k \geq 2n-3,  k \not= n;n_n \geq n &\hbox{when} \hskip .3cm p=n^2-2n+4\\[.2pc]
n_k \geq n^2+1-p,  &\hbox{when} \hskip .3cm n^2-2n+5 \\[.2pc]
~~k \not =n;n_n \geq n &\hskip 1cm \leq p \leq n^2-n+1\\[.2pc]
n_k \geq n^2+1-p, &\hbox{when} \hskip .3cm n^2-n+2 \leq p \leq n^2\\[.2pc]
 ~~ k \not = n;n_n \geq n^2+1-p
\end{cases}.
\end{equation*}
\item Let $G$ be a simply connected algebraic group of type $D_n$. Let
$\lambda = \sum n_k \omega_k.$ Then
\begin{equation*}
\hskip -.55cm H^{p,q}(G/B, \L({\lambda})) =0, \forall ~q\geq1,
\end{equation*}
if
\begin{equation*}
\hskip -.55cm \begin{cases}
n_k \geq p &\hbox{when} \hskip .3cm 0\leq p \leq 2n-4\\[.2pc]
n_k \geq 2n-3 &\hbox{when} \hskip .3cm 2n-3 \leq p \leq n^2-3n+4\\[.2pc]
n_k \geq n^2-n+1-p &\hbox{when} \hskip .3cm
n^2-3n+5 \leq p \leq n^2-n
\end{cases}.
\end{equation*}
\end{enumerate}

\begin{enumerate}
\setcounter{enumi}{4}
\renewcommand\labelenumi{{\rm (\Alph{enumi}$_{6}$)}}
\leftskip .55pc

\item Let $G$ be a simply connected algebraic group of type $E_6$.
Let $\lambda = \sum n_k \omega_k.$ Then
\begin{equation*}
\hskip -.55cm H^{p,q}(G/B, \mathcal{ L}({\lambda})) =0, \forall ~q\geq 1,
\end{equation*}
if
\begin{equation*}
\hskip -.55cm \begin{cases}
n_k \geq p  &\hbox{when} \hskip .3cm 0\leq p \leq 10\\[.2pc]
n_k \geq 11 &\hbox{when} \hskip .3cm 11 \leq p \leq 26\\[.2pc]
n_k \geq 37-p &\hbox{when} \hskip .3cm 27 \leq p \leq 36
\end{cases}.
\end{equation*}
\end{enumerate}

\begin{enumerate}
\setcounter{enumi}{4}
\renewcommand\labelenumi{{\rm (\Alph{enumi}$_{7}$)}}
\leftskip .55pc

\item Let $G$ be a simply connected algebraic group of type $E_7$. Let
$\lambda = \sum n_k \omega_k.$ Then
\begin{equation*}
\hskip -.55cm H^{p,q}(G/B, \L({\lambda})) =0, \forall ~q\geq 1,
\end{equation*}
if
\begin{equation*}
\hskip -.55cm \begin{cases}
      n_k \geq p  &\hbox{when} \hskip .3cm 0\leq p \leq 16\\[.2pc]
      n_k \geq 17  &\hbox{when} \hskip .3cm 17 \leq p \leq 47\\[.2pc]
      n_k \geq 64-p &\hbox{when} \hskip .3cm 48 \leq p \leq 63
\end{cases}.
\end{equation*}
\end{enumerate}

\begin{enumerate}
\setcounter{enumi}{4}
\renewcommand\labelenumi{{\rm (\Alph{enumi}$_{8}$)}}
\leftskip .55pc

\item Let $G$ be a simply connected algebraic group of type $E_8$.
Let $\lambda = \sum n_k \omega_k.$ Then
\begin{equation*}
\hskip -.55cm H^{p,q}(G/B, \L({\lambda})) =0, \forall ~q\geq 1,
\end{equation*}
if
\begin{equation*}
\hskip -.55cm \begin{cases}
      n_k \geq p  &\hbox{when} \hskip .3cm 0\leq p \leq 28\\[.3pc]
      n_k \geq 29  &\hbox{when} \hskip .3cm 29 \leq p \leq 92\\[.3pc]
      n_k \geq 121-p &\hbox{when} \hskip .3cm 93 \leq p \leq 120
\end{cases}.
\end{equation*}
\end{enumerate}

\begin{enumerate}
\setcounter{enumi}{5}
\renewcommand\labelenumi{{\rm (\Alph{enumi}$_{4}$)}}
\leftskip .55pc
\item Let $G$ be a simply connected algebraic group of type $F_4$. Let
$\lambda =\sum n_k \omega_k.$ Then
\begin{equation*}
\hskip -.55cm H^{p,q}(G/B, \L({\lambda})) =0, \forall ~q\geq 1,
\end{equation*}
if
\begin{equation*}
\hskip -.55cm \begin{cases}
n_1,n_2 \geq p; n_3,n_4 \geq 2p-1 &\hbox{when} \hskip .3cm 0\leq p \leq
4\\[.3pc]
n_1,n_2 \geq p; n_3,n_4 \geq p+3 &\hbox{when} \hskip .3cm 5 \leq p
\leq7\\[.3pc]
n_1,n_2 \geq 8; n_3,n_2\geq 11 = h_k-1 &\hbox{when} \hskip .3cm 8 \leq p
\leq 17\\[.3pc]
n_1,n_2 \geq 25-p; &\hbox{when} \hskip .3cm 18\leq p \leq 20\\[.3pc]
n_3,n_4 \geq 45-2p \\[.3pc]
n_k \geq 25-p &\hbox{when} \hskip .3cm 21 \leq p \leq 24
\end{cases}.
\end{equation*}
\end{enumerate}

\begin{enumerate}
\setcounter{enumi}{6}
\renewcommand\labelenumi{{\rm (\Alph{enumi}$_{2}$)}}
\leftskip .55pc
\item Let $G$ be a simply connected algebraic group of type $G_2$. Let
$\lambda = n \omega_1+m \omega_2$. Then
\begin{equation*}
\hskip -.55cm H^{p,q}(G/B, \L({\lambda})) = 0,\forall ~q\geq 1,
\end{equation*}
if
\begin{equation*}
\hskip -.55cm \begin{cases}
      n,m \geq 1 &\hbox{when} \hskip .3cm p =1,6\\[.3pc]
      n\geq 2, m\geq 4 &\hbox{when} \hskip .3cm p=2,5 \\[.3pc]
      n\geq 3, m\geq 5 &\hbox{when} \hskip .3cm p=3,4
\end{cases}.
\end{equation*}
\end{enumerate}
\end{propo}

\begin{proof}
We consider the set of weights of $\Lambda ^p \frak{n^-}\otimes
k_{\lambda}$. By our notation this set is $\Phi ^-_p +\{\lambda\}$. Our
aim is to make all $\nu \in \Phi ^-_p +\{\lambda\}$ either dominant or
$\lambda+\nu$ singular for suitable choice of $\lambda$. For that we
give condition on $(\lambda, \alpha^{\vee}_i)$ so that $(\nu,
\alpha_i^{\vee})\geq-1$ for all $ \alpha_i \in \Delta$. To avoid
complications we include only simple singularities, i.e. we give
conditions on $(\lambda, \alpha_i^{\vee})$ so that either $\nu \in
X(T)^+$, or there exists a simple root $\alpha$ such that $(\nu+\rho,
\alpha^{\vee}) =0.$

If we know in each column of the `positive roots matrix', what the
positive entries, the negative entries, and the zeros are, and how many
times they appear, then we can compute easily what the maximum positive
value is which can occur in $(\mu, \alpha_i^{\vee})$ for all
$\alpha_i^{\vee}\in \Delta$, where $\mu $ runs over $\Phi ^+_p$. Suppose
that this number is $n_i$, then we demand that $\lambda$ satisfy the
inequalities $(\lambda,\alpha_i^{\vee}) \geq n_i -1 \forall i$.

The following observation will be useful in writing down the conditions
stated in the proposition. In the positive root matrix, suppose there
are $m_1$ positive entries and $m_2$ negative entries and $m_3$ zeros in
the $i$th column. Then the maximum among the positive entries in the
$i$th column of the matrix corresponding to $\Phi^+_p$ increases when $p$
increases from zero to $m_1$ (how quickly it increases depends on
what the positive entries are (with repetition)); it will be
non-decreasing when $p$ increases from $m_1$ to $m_1+m_3$, and then it
starts decreasing (again, how quickly it decreases depends on how many
negative entries there are (with repetition)).

We give the needed information about the positive roots matrix in the
following table.\vspace{1pc}

%\begin{table}[h]
%\processtable{}
{\begin{tabular}{@{}l@{\hskip .5cm}c@{\hskip .5cm}|c@{\hskip .5cm}c@{\hskip .5cm}c@{\hskip .5cm}c@{\hskip .5cm}c@{\hskip .5cm}c@{\hskip .5cm}c@{}}\cline{1-9}
 & & & & & & & &\\[-.7pc]
 & &0 &1 &2 &3 &$-1$ &$-2$ &$-3$\\
 & & & & & & & &\\[-1.1pc]
\multicolumn{9}{c@{}@{}}{\hskip -.5pc\rule[.4pc]{11.2cm}{.009cm}}\\[-.4pc]
 & & & & & & & &\\[-.7pc]
$A_{n}$ &   &$\displaystyle\frac{n^2-3n+2}{2}$ &$n-1$ &1 &0 &$n-1$ &0 &0\\[.8pc]
$B_n$   &Short  &$n^2-2n+1$ &0 &$n$ &0 &0 &$n-1$ &0\\[.5pc]
 &Long &$n^2-4n+5$ &$2n-3$ &1 &0 &$2n-3$ &0 &0\\[.5pc]
$C_n$ &Short &$n^2-4n+5$ &$2n-4$ &2 &0 &$2n-4$ &1 &0 \\[.5pc]
&Long &$n^2-2n+1$ &$n-1$ &1 &0 &$n-1$ &0 &0\\[.5pc]
$D_n$ & &$n^2-5n+7$ &$2n-4$ &1 &0 &$2n-4$ &0 &0\\[.5pc]
$E_6$ & &15 &10 &1 &0 &10 &0 &0\\[.5pc]
$E_7$ & &30 &16 &1 &0 &16 &0 &0\\[.5pc]
$E_8$ & &63 &28 &1 &0 &28 &0 &0\\[.5pc]
$F_4$ &Short &9 &4 &4 &0 &4 &3 &0\\[.5pc]
      &Long &9 &7 &1 &0 &7 &0 &0\\[.5pc]
$G_2$ &Short &1 &1 &1 &1 &1 &0 &1\\[.5pc]
 &Long &1 &2 &1 &0 &2 &0 &0\\[.15pc]\hline
\end{tabular}}{}
%\end{table}

\hfill $\BBox\,$\vspace{.1pc}
\end{proof}

\setcounter{rem}{2}
\begin{rem}
{\rm We do not present all the computations of the table mainly for want of
space. However, we give some computations for classical groups in the
next section. We give explicitly the matrix corresponding to $\Phi^+$
for all the exceptional groups in Appendix A.}
\end{rem}

\begin{rem}
{\rm Let $\lambda= \sum n_k \omega_k$. We gave conditions on $n_{k}$'s
depending on $p$ so as to have $H^{p,q}(G/B, \L (\lambda))=0, \forall q
\geq 1$. But conversely the vanishing $H^{p,q}(G/B, \L (\lambda))=0,
\forall q \geq 1$ need not imply that the $n_{k}$'s satisfy such
conditions. For example let $G$ be of type $E_6$ and let
$p=\dim(G/B)-1.$ If $\lambda = 2\rho-\omega_2$ or $\lambda = \rho$ then
$\lambda$ will not satisfy the conditions given in the theorem but
$H^{p,q}(G/B, \L (\lambda))=0, \forall q \geq 1$.}
\end{rem}

\setcounter{coro}{4}
\begin{coro}$\left.\right.$\vspace{.5pc}

\noindent $H^{p,q}(G/B, \L(\lambda))=0, \forall ~q\geq 1, \forall
p\geq0,$ when $n_k \geq h_{\alpha_k}-1, ~ \forall k.$
\end{coro}

\begin{proof}
It follows from Proposition 2 and Theorem 1. \hfill $\BBox\,$
\end{proof}

The following corollary is a weaker version of the above corollary.

\begin{coro}$\left.\right.$\vspace{.5pc}

\noindent $H^{p,q}(G/B, \L(\lambda))=0, \forall ~q\geq 1, \forall
p\geq0,$ when $n_k \geq h-1, ~ \forall k.$
\end{coro}

\setcounter{rem}{6}
\begin{rem}
{\rm We could directly prove Corollary 6 along the same lines of the
proof of Theorem~1 with the following observations. This is due to
Andersen and Jantzen \cite{ander.jant}. Let $\nu$ be a
weight of $\Lambda^j\frak{n}^-,$ for some $j$. Then $\max _{\alpha \in
\Phi}|(\nu+\rho, \alpha^{\vee})|\leq h-1$.}
\end{rem}

\begin{proof}
Recall that $\rho$ is the half sum of positive roots. Since $\nu$ is the
sum of some distinct negative roots, we have $\nu +\rho =
\frac{1}{2}\sum _{\alpha \in S}\alpha -\frac{1}{2}\sum _{\alpha \in
S^{\prime}}\alpha$ where $S^{\prime} \subset \Phi ^+$ such that $\nu =
\sum _{{\small \alpha \in S^{\prime}}} -\alpha,$ and $S = \Phi^+ -
S^{\prime}$. Since $W$ acts on $\Phi$, $w(\nu +\rho)$ will be of the
above form and hence $w(\nu +\rho)\leq \rho,~\forall ~ w \in W.$ Now fix
$w$ such that $w(\nu +\rho)\in X(T)^+. $

One has
\begin{align*}
\mathop{\max}\limits_{\alpha \in \Phi}|(\nu+\rho, \alpha^{\vee})| &=
\mathop{\max}\limits_{\alpha \in \Phi}|(w(\nu+\rho), \alpha^{\vee})|\\[.2pc]
&=\mathop{\max}\limits_{\alpha \in \Phi}(w(\nu+\rho), \alpha^{\vee})\leq \
\mathop{\max}\limits_{\alpha \in \Phi}(\rho, \alpha^{\vee})= h-1
\end{align*}
which gives the inequality quoted in the above remark.\hfill $\BBox\,$\vspace{-.5pc}
\end{proof}

\section{A remark on the Coxeter number}

From the above table, we can easily obtain the following lemma.%\vspace{-.3pc}

\setcounter{lem}{7}
\begin{lem}$\left.\right.$
\begin{enumerate}
\renewcommand\labelenumi{\rm (\alph{enumi})}
\leftskip .1pc
\item If $\alpha_i${\rm ,} $\alpha_j \in \Delta$ such that
$|\alpha_i|=|\alpha_j|$ $\hbox{then}\ h_{\alpha_i}=h_{\alpha_j}$.

\item If $\alpha_i${\rm ,} $\alpha_j \in \Delta$ such that
$|\alpha_i|\leq |\alpha_j|$ $\hbox{then}\ h_{\alpha_i} \geq
h_{\alpha_j}.$\vspace{-.8pc}
\end{enumerate}
\end{lem}

\begin{proof}
As we mentioned, the proof can be read off from the table. However, we
give a conceptual proof for (a). Let $\alpha \in \Phi^+$. For any
$\gamma \in \Phi^+ -\{\alpha\},$ consider the $\alpha$-string
$\gamma-p\alpha, \ldots , \gamma+q\alpha$ through $\gamma$, where $p-q =
(\gamma, \alpha^{\vee})$. It follows that
$(\gamma-(p-i)\alpha,\alpha^{\vee})=-p-q+2i=-(\gamma+(q-i)\alpha,
\alpha^{\vee})$. Since $\gamma \in \Phi^+ -\{\alpha\},$ is in a unique
$\alpha$-string of maximum length, we have $h_{\alpha} = \sum _{\gamma
\in \Phi^+_{\alpha}}(\gamma, \alpha^{\vee}) = \frac{1}{2}(2+\sum
_{\gamma \in \Phi^+}|(\gamma, \alpha^{\vee})|)=\frac{1}{4}(4+\sum
_{\gamma \in \Phi}|(\gamma, \alpha^{\vee})|).$ Now since the pairing
$(,)$ is $W$-invariant and since $W$ acts transitively on the set of
roots of $G$ of the same length, (a) follows. \hfill $\BBox\,$
\end{proof}

\setcounter{theorr}{8}
\begin{therr}[\!]
Let $\overline h$ be $\max_{\alpha \in \Delta}\{h_{\alpha}\}$. The
number $\overline h$ coincides with the Coxeter number $h$.
\end{therr}

\begin{proof}
The proof we give uses the Dynkin classification. The proof follows from
the table. Since we have not given computation of the table, here we
give the proof for certain cases in classical type. Appendix~A gives
complete proofs for all the exceptional groups.

\begin{enumerate}
\renewcommand{\labelenumi}{\arabic{enumi}.}
\leftskip -.2pc
\item $A_n$ type

Any root $\gamma \in \Phi^+$ is a linear combination of simple roots
with coefficient 0 or 1. Let $\{\alpha_1,~ \alpha_2 , \ldots ,~
\alpha_n \}$ be a set of simple roots (with the usual order). We write
the positive roots as an ordered linear combination of simple roots
(i.e. if $\alpha_1$ appears in the linear combination then it should
appear first and then $\alpha_2$ and so on).
Let $\alpha_i$ be a simple root. It is clear that for $\gamma \in
\Phi^+, ~ (\gamma, \alpha_i^{\vee}) > 0 $ (i.e. when we write $\gamma$
in terms of fundamental weights, the coefficient of $\omega_i$ is
strictly positive) iff $\gamma$ begins with $\alpha_i$ or ends with
$\alpha_i$. The number of positive roots beginning with $\alpha_i$ is $
n+1-i$ (the roots are $\alpha_i, ~ \alpha_i+\alpha_{i+1}, \ldots ,
\alpha_i+\alpha_{i+1}+ \cdots + \alpha_n)$, and the number of positive
roots ending with $\alpha_i$ is $i$ (the roots are $\alpha_i,~
\alpha_{i-1}+\alpha_i, \ldots , ~ \alpha_1+\alpha_2+ \cdots +
\alpha_i)$. Note that we have counted $\alpha_i$ twice. Therefore the
total number of positive roots with $\omega_i$th coefficient strictly
positive is $n$. Among all such positive roots, only $\alpha_i$ will have
$\omega_i$th coefficient to be $2$, all other positive roots will have
coefficient to be $1$. So $h_{\alpha_i}=n+1, ~\forall i$, hence
$\overline h=n+1 =h$.\vspace{.5pc}

\item $B_n$ type

We compute $h_{\alpha_n}$. The positive roots with $\omega_n$th
coefficient strictly positive are $\alpha_n, ~ \alpha_{n-1}+2\alpha_n,
~\alpha_{n-2}+\alpha_{n-1}+2\alpha_n, \ldots ,~ \alpha_1 + \cdots
+\alpha_{n-1}+2\alpha_n,$ and they all will have $\omega_n$th
coefficient $2$, hence $\overline h = h_{\alpha_n}=2n=h$.\vspace{.5pc}

\item $C_n$ type

We compute $h_{\alpha_1}.$ The positive roots with $\omega_1$th
coefficient strictly positive are $\alpha_1, ~
\alpha_1+\alpha_2,~\alpha_1+\alpha_2+\alpha_3, \ldots ,~ \alpha_1 +
\cdots+ \alpha_n$ (total number is $n$), and
$\alpha_1+\alpha_2+2\alpha_3+ \cdots +2\alpha_{n-1}+\alpha_n, ~
\alpha_1+\alpha_2+\alpha_3+2\alpha_4 +\cdots+ 2\alpha_{n-1} + \alpha_n,
 \ldots,~\alpha_1+\alpha_2+\cdots +\alpha_{n-2}+ 2\alpha_{n-1}+
\alpha_n$ (total number is $n-3$, and $2\alpha_1+2\alpha_2+\cdots
+2\alpha_{n-1}+ \alpha_n.$ All of them will have $\omega_1$th
coefficient $1$, except $\alpha_1$ and $2\alpha_1+2\alpha_2+\cdots
+2\alpha_{n-1}+ \alpha_n,$ and $\omega_1$th coefficient of these two
weights will be $2,$ hence $\overline h = h_{\alpha_n}=2n=h$.\vspace{.5pc}

\item $D_n$ type

We compute $h_{\alpha_1}$.  The positive roots with $\omega_1$th
coefficient strictly positive are $\alpha_1, ~ \alpha_1+\alpha_2, ~
\alpha_1+\alpha_2+ \cdots +\alpha_n$ (total number is $n$),
and $\alpha_1+\alpha_2+  \cdots + \alpha_{n-2}+\alpha_n$, and
$\alpha_1+\alpha_2+2\alpha_3+\cdots+2\alpha_{n-2}+\alpha_{n-1}+\alpha_n,~
\alpha_1+\alpha_2+\alpha_3+2\alpha_4 + \cdots+2\alpha_{n-2}+\alpha_{n-1}+
\alpha_n,\ldots,~ \alpha_1+\alpha_2+\cdots+
\alpha_{n-3}+2\alpha_{n-2}+\alpha_{n-1}+\alpha_n$ (total number is
$n-4$).  All of them will have $\omega_1$th coefficient $1$ except
$\alpha_1,$ and $\alpha_1$ will have $2$, and hence
$h_{\alpha_1}=\overline{h}=2n-2=h$. \hfill $\BBox\ \,$\vspace{-1pc}
\end{enumerate}
\end{proof}

\section{Remarks on Bott vanishing property}

In this section we prove that the Bott vanishing property fails to hold
on generalised flag varieties. More precisely, we prove that for a
simple, simply connected algebraic group $G$ of rank $\geq 2$, there
exists ample line bundle $L_{\lambda}$ with $H^{d-1,1}(G/B,
\L({\lambda}))\not = 0$, where $d=\dim(G/B)$. Since $G/B$ is
defined over ${\Bbb Z},$ the non-vanishing results hold even in prime
characteristic for large primes.

\setcounter{lem}{9}
\begin{lem}
Let $V$ be a finite dimensional $B$-module with a one-dimensional
$B$-submodule $k_{\lambda}$ and let $U$ be the quotient module
$V/k_{\lambda}.$ If $\lambda \in X(T)-W\cdot X(T)^+ $ then $H^i(G/B,
\L(V))\simeq H^i(G/B, \L(U)), \forall i$.
\end{lem}

\begin{proof}
Consider the short exact sequence of $B$-modules
$0\rightarrow K_{\lambda} \rightarrow V \rightarrow U
\rightarrow0$.  This induces a long exact sequence of cohomology
modules
\begin{align*}
0 &\rightarrow H^0(G/B, \L(\lambda)) \rightarrow H^0(G/B,
\L(V))\rightarrow H^0(G/B, \L(U))\\[.2pc]
&\rightarrow H^1(G/B, \L(\lambda))\rightarrow H^1(G/B, \L(V))
\rightarrow H^1(G/B,  \L(U))\\[.2pc]
&\rightarrow H^2(G/B, \L(\lambda))
\rightarrow H^2(G/B, \L(V))\rightarrow H^2(G/B, \L(U))\rightarrow \cdots
\end{align*}
By the Borel--Weil--Bott theorem we know that $H^i(G/B, \L(\lambda))=0$
for all $i$.\hfill $\BBox\,$
\end{proof}

\begin{lem}$\left.\right.$
\begin{enumerate}
\renewcommand\labelenumi{\rm (\roman{enumi})}
\leftskip .3pc
\item Let $G$ be simple group with {\rm rank}$(G) \geq 4${\rm ,} other than type
$F_4$ and $D_n$. Let $\lambda=2\rho- \alpha_{n-2}-\alpha_{n-1}$. Then
$\lambda$ is strictly dominant and the weights of
$\Lambda^{d-1}\frak{n}^-\otimes k_{\lambda}$ all lie in $X(T)-W\cdot
X(T)^+$ except $-\alpha_{n-1}, -\alpha_{n-2}, 0$. More precisely{\rm ,} if
$\mu$ is a weight of $\Lambda^{d-1}\frak{n}^-\otimes k_{\lambda}$ and
$\mu \not \in \{ -\alpha_{n-2}, -\alpha_{n-1}, 0\}${\rm ,} then there exists
a $\nu \in \{ \alpha_{n-2}, \alpha_{n-1}, \alpha_{n-2}+ \alpha_{n-1}
\}$ such that $(\mu +\rho, \nu^{\vee})=0$.

\item Let $G$ be a simple group of {\rm rank} $2$ or {\rm rank} $3$ or of type
$F_4$. Let $\lambda=2\rho- \alpha_1-\alpha_2$. Then $\lambda$ is
strictly dominant and the weights of $\Lambda^{d-1}\frak{n}^-\otimes
k_{\lambda}$ all lie in $X(T)-W\cdot X(T)^+$ except $-\alpha_2,
-\alpha_1, 0$ {\rm (}for $G_2${\rm ,} in addition to this we will have $(0,1)$
also as a non-singular weight{\rm )}. More precisely{\rm ,} if $\mu$ is a weight of
$\Lambda^{d-1}\frak{n}^-\otimes k_{\lambda}$ and $\mu \not \in \{
-\alpha_1, -\alpha_2, 0\}${\rm ,} then there exists a $\nu \in \{ \alpha_1,
\alpha_2, \alpha_1+ \alpha_2 \}$ such that $(\mu +\rho,
\nu^{\vee})=0$.

\item Let $G$ be $D_n$ type. Let $\lambda = 2 \rho -
\alpha_{n-3} - \alpha_{n-2}$. Then $\lambda$ is strictly dominant and the
weights of $\Lambda^{d-1}\frak{n}^- \otimes k_{\lambda}$ all lie in
$X(T)-W\cdot X(T)^+$ except $-\alpha_{n-2}, -\alpha_{n-3}, 0$. More
precisely{\rm ,} if $\mu$ is a weight of $\Lambda^{d-1}\frak{n}^-\otimes
k_{\lambda}$ and $\mu \not \in \{ -\alpha_{n-3}, -\alpha_{n-2}, 0\},$
then there exists a $\nu \in \{ \alpha_{n-3}, \alpha_{n-2},
\alpha_{n-3}+ \alpha_{n-2} \}$ such that $(\mu +\rho, \nu^{\vee})=0$.
\end{enumerate}
\end{lem}

\begin{proof}
For the exceptional groups, we could verify the lemma directly from
the tables given in the Appendix.

For the classical groups we essentially need to prove the statements for
rank $4$ and rank $5$ groups for $D_n$ type, for other types rank $\leq
4$ groups. For these groups it is very easy to check it by hand. For the
convenience of the reader we give the table of positive roots in these
cases. Then we can easily check the statements.

{\abb
$$\begin{tabular}{@{}l@{\,\ }c@{\ \ \ \ \ \,}r@{\ \ \,}r@{\ \ \,}r@{\ \ \,}r l@{\,\ }c@{\ \ \ \ \ \,}r@{\ \ \,}r@{\ \ \,}r@{\ \ \,}r l@{\,\ }c@{\ \ \ \ \ \,}r@{\ \ \,}r@{\ \ \,}r@{\ \ \,}r@{}}
\multicolumn{6}{c}{$A_4$} &\multicolumn{6}{c}{$B_4$} &\multicolumn{6}{c}{$C_4$}\\
&&&&&&&&&&&\\
(1 0 0 0)&\myfavouritearrow &($\phantom{-}$2 & $-$1 & 0 & 0) &
(1 0 0 0)&\myfavouritearrow &($\phantom{-}$2 & $-$1 & 0 & 0) &
(1 0 0 0)&\myfavouritearrow &($\phantom{-}$2 & $-$1 & 0 & 0) \\

(0 1 0 0)&\myfavouritearrow &($-$1 & 2 & $-$1 & 0) &
(0 1 0 0)&\myfavouritearrow &($-$1 & 2 & $-$1 & 0) &
(0 1 0 0)&\myfavouritearrow &($-$1 & 2 & $-$1 & 0) \\

(0 0 1 0)&\myfavouritearrow &($\phantom{-}$0 & $-$1 & 2 & $-$1) &
(0 0 1 0)&\myfavouritearrow &($\phantom{-}$0 & $-$1 & 2 & $-$2) &
(0 0 1 0)&\myfavouritearrow &($\phantom{-}$0 & $-$1 & 2 & $-$1) \\

(0 0 0 1)&\myfavouritearrow &($\phantom{-}$0 & 0 & $-$1& 2) &
(0 0 0 1)&\myfavouritearrow &($\phantom{-}$0 & 0 & $-$1& 2) &
(0 0 0 1)&\myfavouritearrow &($\phantom{-}$0 & 0 & $-$2& 2) \\

(1 1 0 0)&\myfavouritearrow &($\phantom{-}$1 & 1 & $-$1 & 0) &
(1 1 0 0)&\myfavouritearrow &($\phantom{-}$1 & 1 & $-$1 & 0) &
(1 1 0 0)&\myfavouritearrow &($\phantom{-}$1 & 1 & $-$1 & 0) \\

(0 1 1 0)&\myfavouritearrow &($-$1 & 1 & 1 &  $-$1)&
(0 1 1 0)&\myfavouritearrow &($-$1 & 1 & 1 &  $-$2) &
(0 1 1 0)&\myfavouritearrow &($-$1 & 1 & 1 &  $-$1) \\

(0 0 1 1)&\myfavouritearrow &($\phantom{-}$0 & $-$1 & 1 & 1) &
(0 0 1 1)&\myfavouritearrow &($\phantom{-}$0 & $-$1 & 1 & 0) &
(0 0 1 1)&\myfavouritearrow &($\phantom{-}$0 & $-$1 & 0 & 1) \\

(1 1 1 0)&\myfavouritearrow &($\phantom{-}$1 & 0 & 1 & $-$1) &
(1 1 1 0)&\myfavouritearrow &($\phantom{-}$1 & 0 & 1 & $-$2) &
(1 1 1 0)&\myfavouritearrow &($\phantom{-}$1 & 0 & 1 & $-$1) \\

(0 1 1 1)&\myfavouritearrow &($-$1 & 1 & 0 & 1)&
(0 1 1 1)&\myfavouritearrow &($-$1 & 1 & 0 & 0)&
(0 1 1 1)&\myfavouritearrow &($-$1 & 1 & $-$1 & 1) \\

(1 1 1 1)&\myfavouritearrow &($\phantom{-}$1 & 0 & 0 & 1) &
(0 0 1 2)&\myfavouritearrow &($\phantom{-}$0 & $-$1 & 0 & 2) &
(0 0 2 1)&\myfavouritearrow &($\phantom{-}$0 & $-$2 & 2 & 0) \\

&&&&&& (1 1 1 1) &\myfavouritearrow &($\phantom{-}$1 & 0 & 0 & 0) &
(1 1 1 1) &\myfavouritearrow &($\phantom{-}$1 & 0 & $-$1 & 1) \\

&&&&&& (0 1 1 2) &\myfavouritearrow &($-$1 & 1 & $-$1 & 2)&
(0 1 2 1) &\myfavouritearrow &($-$1 & 0 & 1 & 0)\\

&&&&&& (1 1 1 2) &\myfavouritearrow &($\phantom{-}$1 & 0 & $-$1 & 2) &
 (1 1 2 1) &\myfavouritearrow &($\phantom{-}$1 & $-$1 & 1 & 0) \\

&&&&&& (0 1 2 2) &\myfavouritearrow &($-$1 & 0 & 1 & 0)&
(0 2 2 1) &\myfavouritearrow &($-$2 & 2 & 0 & 0)\\

&&&&&& (1 1 2 2) &\myfavouritearrow &($\phantom{-}$1 & $-$1 & 1 & 0) &
 (1 2 2 1) &\myfavouritearrow &($\phantom{-}$0 & 1 & 0 & 0) \\

&&&&&& (1 2 2 2) &\myfavouritearrow &($\phantom{-}$0 & 1 & 0 & 0) &
 (2 2 2 1) &\myfavouritearrow &($\phantom{-}$2 & 0 & 0 & 0) \\
\end{tabular} $$

$$\begin{tabular}{@{}l@{\ \,}c@{\ \ \ \ \ \,}r@{\ \ \,}r l@{\ \,}c@{\ \ \ \ \ \,}r@{\ \ \,}r@{\ \ \,}r l@{\ \,}c@{\ \ \ \ \ \,}r@{\ \ \,}r@{\ \ \,}r l@{\ \,}c@{\ \ \ \ \ \,}r@{\ \ \,}r@{\ \ \,}r@{}}
\multicolumn{4}{c}{${\boldmath B_2}$} &\multicolumn{5}{c}{${\boldmath A_3}$}
&\multicolumn{5}{c}{${\boldmath B_3}$} &\multicolumn{5}{c}{${\boldmath C_3}$}\\
& & & & & & & & & & & & &\\
(1 0) &\myfavouritearrow &($\phantom{-}$2 & $-$2) &
(1 0 0) &\myfavouritearrow &($\phantom{-}$2 & $-$1 &0) &
(1 0 0) &\myfavouritearrow &($\phantom{-}$2 &$-$1 &0)&
(1 0 0) &\myfavouritearrow &($\phantom{-}$2 &$-$1 &0) \\

(0 1) &\myfavouritearrow &($-$1 & 2)&
(0 1 0) &\myfavouritearrow &($-$1 &2& $-$1) &
(0 1 0) &\myfavouritearrow &($-$1 & 2& $-$2)&
(0 1 0) &\myfavouritearrow &($-$1 & 2& $-$1) \\

(1 1) &\myfavouritearrow &($\phantom{-}$1 & 0)&
(0 0 1) &\myfavouritearrow &($\phantom{-}$0 & $-$1& 2) &
(0 0 1) &\myfavouritearrow &($\phantom{-}$0 & $-$1 & 2)&
(0 0 1) &\myfavouritearrow &($\phantom{-}$0 &$-$2 &  2) \\

(1 2) &\myfavouritearrow &($\phantom{-}$0 & 2)&
(1 1 0) &\myfavouritearrow &($\phantom{-}$1 &1 & $-$1) &
(1 1 0) &\myfavouritearrow &($\phantom{-}$1 &1 &$-$2)&
(1 1 0) &\myfavouritearrow &($\phantom{-}$1& 1 & $-$1) \\

& & & & (0 1 1) &\myfavouritearrow &($-$1 &1 & 1) &
(0 1 1) &\myfavouritearrow &($-$1 & 1& 0)&
(0 1 1) &\myfavouritearrow &($-$1 &  0 & 1) \\

& & & & (1 1 1) &\myfavouritearrow &($\phantom{-}$1 &0 & 1) &
(1 1 1) &\myfavouritearrow &($\phantom{-}$1 & 0 & 0) &
(1 1 1) &\myfavouritearrow &($\phantom{-}$1 & $-$1 &1) \\

&&& & & &  & &  &
(0 1 2) &\myfavouritearrow &($-$1 & 0& 2)&
(0 2 1) &\myfavouritearrow &($-$2& 2 &  0) \\

&&& & & &  & &  &
(1 1 2) &\myfavouritearrow &($\phantom{-}$1 & $-$1& 2)&
(1 2 1) &\myfavouritearrow &($\phantom{-}$0& 1 &  0) \\

&&& & & &  & &  &
(1 2 2) &\myfavouritearrow &($\phantom{-}$0 & 1& 0)&
(2 2 1) &\myfavouritearrow &($\phantom{-}$2& 0 &  0) %\\
\end{tabular}$$
$$\begin{tabular}{@{}lc@{\ \ \ \ \ \ \ \ \ \ \,}rrrr r lc@{\ \ \ \ \ \ \ \ \ \ \,}rrrrr@{}}
%\begin{tabular}{rrrrrrrrrrrrrr}
& &  ${\boldmath D_4}$ & & & & & & & & ${\boldmath D_5}$ & & & \\
& &  & & & & &  & & & \\[-.4pc]
(1 0 0 0)& \myfavouritearrow & ($\phantom{-}$2 & $-$1 &0& 0)& &
(1 0 0 0 0) & \myfavouritearrow & ($\phantom{-}$2 & $-$1& 0 & 0 &0) \\
(0 1 0 0)& \myfavouritearrow & ($-$1 & 2 & $-$1 &$-$1)& &
(0 1 0 0 0) & \myfavouritearrow & ($-$1& 2 & $-$1& 0 & 0) \\
(0 0 1 0)& \myfavouritearrow & ($\phantom{-}$0 & $-$1 &2& 0)& &
(0 0 1 0 0) & \myfavouritearrow & ($\phantom{-}$0 & $-$1& 2 & $-$1 &$-$1) \\
(0 0 0 1)& \myfavouritearrow & ($\phantom{-}$0 & $-$1 &0& 2)& &
(0 0 0 1 0) & \myfavouritearrow & ($\phantom{-}$0 & 0&$-$1& 2 & 0) \\
(1 1 0 0)& \myfavouritearrow & ($\phantom{-}$1 & 1& $-$1 &$-$1)& &
(0 0 0 0 1) & \myfavouritearrow & ($\phantom{-}$0 & 0& $-$1& 0 & 2) \\
(0 1 1 0)& \myfavouritearrow & ($-$1 & 1 & 1 &$-$1)& &
(1 1 0 0 0) & \myfavouritearrow & ($\phantom{-}$1&1 & $-$1& 0 & 0) \\
(0 1 0 1)& \myfavouritearrow & ($-$1 & 1 & $-$1 & 1)& &
(0 1 1 0 0) & \myfavouritearrow & ($-$1& 1&1& $-$1&$-$1) \\
(1 1 1 0)& \myfavouritearrow & ($\phantom{-}$1 & 0& 1 &$-$1)& &
(0 0 1 1 0) & \myfavouritearrow & ($\phantom{-}$0 & $-$1& 1 & 1 &$-$1) \\
(1 1 0 1)& \myfavouritearrow & ($\phantom{-}$1 & 0& $-$1 &1)& &
(0 0 1 0 1) & \myfavouritearrow & ($\phantom{-}$0 & $-$1& 1 & $-$1 &1) \\
(0 1 1 1)& \myfavouritearrow & ($-$1 & 0 & 1 & 1)& &
(1 1 1 0 0) & \myfavouritearrow & ($\phantom{-}$1 & 0& 1 & $-$1&$-$1) \\
(1 1 1 1)& \myfavouritearrow & ($\phantom{-}$1 & $-$1& 1 &1)& &
(0 1 1 1 0) & \myfavouritearrow & ($-$1 & 1& 0 & 1 &$-$1) \\
(1 2 1 1)& \myfavouritearrow & ($\phantom{-}$0 & 1& 0 &0)& &
(0 0 1 1 1) & \myfavouritearrow & ($\phantom{-}$0 & $-$1& 0 & 1 &1) \\
&&&&&&& (0 1 1 0 1)& \myfavouritearrow & ($-$1& 1 & 0 & $-$1 & 1)\\
&&&&&&& (1 1 1 1 0)& \myfavouritearrow &($\phantom{-}$1& 0& 0& 1&$-$1)\\
&&&&&&& (0 1 1 1 1)& \myfavouritearrow & ($-$1 & 1 & $-$1 &1 & 1)\\
&&&&&&& (1 1 1 0 1)& \myfavouritearrow &($\phantom{-}$1& 0 & 0 & $-$1 & 1)\\
&&&&&&& (1 1 1 1 1)& \myfavouritearrow &($\phantom{-}$1& 0&$-$1 & 1&1)\\
&&&&&&& (0 1 2 1 1)& \myfavouritearrow & ($-$1 & 0 & 1 &0 & 0)\\
&&&&&&& (1 1 2 1 1)& \myfavouritearrow &($\phantom{-}$1& $-$1 & 1 & 0 &0)\\
&&&&&&& (1 2 2 1 1)& \myfavouritearrow &($\phantom{-}$0& 1&0 & 0&0)%\\
 \end{tabular} $$}

Let $G$ be a classical group (other than $D_n$ type). We prove the lemma
for the higher rank groups. Similar proof holds for $D_n$ type also. Let
$G$ be of rank $\geq 4$. Note that the weights of
$\Lambda^{d-1}\frak{n}^-\otimes k_{\lambda}$ are $\alpha
-\alpha_{n-2}-\alpha_{n-2}, ~ \alpha \in \Phi^+$. Let $\alpha
=\beta+\gamma$, where $\beta$ involves simple roots $\alpha_1$,
$\ldots$, $\alpha_{n-4}$ and $\gamma$ involves simple roots
$\alpha_{n}$, $\alpha_{n-1}$, $\alpha_{n-2}$, $\alpha_{n-3}$. Let $\nu
\in \{\alpha_{n-2},~ \alpha_{n-1}, ~\alpha_{n-2}+\alpha_{n-1}\}$ then
$(\beta, \nu^{\vee})=0$. Now $(\alpha -\alpha_{n-2}-\alpha_{n-2}+\rho,
\nu ^{\vee})= (\beta, \nu^{\vee}) + (\gamma - \alpha_{n-2} -
\alpha_{n-2} + \rho, \nu^{\vee})$. We could `think of' $\gamma$ as
positive root of rank $\leq 4$ group, which proves the lemma. (In the
case of type $B_n$, for the root $\alpha = \alpha_1 + 2\alpha_2 +
2\alpha_3 + \cdots + 2\alpha_n$, $\gamma$ will be $2\alpha_{n-3} +
2\alpha_{n-2} + 2\alpha_{n-1} + 2\alpha_n$, which will not be a root of
$B_4$, but we can easily see that for $\nu = \alpha_{n-2} +
\alpha_{n-1}$, $(\gamma - \alpha_{n-2} - \alpha_{n-1} + \rho,
\nu^{\vee})$\break $ = 0$.)\hfill $\BBox\,$
\end{proof}

\setcounter{theorr}{11}
\begin{therr}[\!]
Let $G$ be an algebraic group with {\rm rank}$(G) \geq 2$. Let $\lambda$ be
as given in the above lemma. Then $H^{d-1,1}(G/B, \L({\lambda}))\not=0$.
\end{therr}

\begin{proof}
Arrange the weights of $\Lambda^{d-1}\frak{n}^-\otimes k_{\lambda}$ as
$\lambda_1, \ldots , \lambda_d$ such that if $s<t$ then $\lambda_t \not
\leq \lambda_s$; in other words, $\lambda_s$ is the lowest weight in the
set of weights \{$ \lambda_{s+1}, \lambda_{s+2}, \ldots ,\lambda_d$\} for
all $s$. There are several ways to arrange the weights with the said
property. Note that the weights of $\Lambda^{d-1}\frak{n}^-\otimes
k_{\lambda}$ are $\alpha-\beta, ~ \alpha \in \Phi^+$ where
$\lambda=2\rho-\beta$. We arrange the weight as follows:
$\alpha_1-\beta$, $\alpha_2-\beta$, $\ldots$, $\alpha_n-\beta$,
$\alpha_1+\alpha_2-\beta$, $\alpha_2+\alpha_3-\beta $, $\ldots$.

If we choose $\lambda$ as given in Lemma~11, the arrangement of weights
of $\Lambda^{d-1} \frak{n} \otimes k_{\lambda}$ will be of the following
form: $\ldots$, $\mu_1$, $\mu_2$, $\ldots$, $0$, $\ldots$. From Lemma
11, it is clear that if there are weights in the dotted places, they all
will lie in $X(T)-W\cdot X(T)^+$. From the Borel--Weil--Bott theorem
$H^i(G/B, \L(\mu_j)) \neq 0$ only when $i=1$, for both $j=1,2$, and
$H^i(G/B, \L(0))\neq 0$ only when $i=0$ (note that there will be one
extra dominant weight contributing to $H^0(G/B)$ after zero weight in
case of $G_2$, the same technique of the proof will do even in this
case). $(\star)$

Let $V_1=(\Lambda^p \frak{n}^-\otimes k_{\lambda})/k_{\lambda_1}$. Then
we have short exact sequence $0\rightarrow k_{\lambda_1} \rightarrow
\Lambda^p \frak{n}^-\otimes k_{\lambda} \rightarrow V_1 \rightarrow0$.
This induces the following long exact sequence of cohomology modules.
\begin{align*}
0 &\rightarrow H^0(G/B,\L({\lambda_1})) \rightarrow
H^{p,0}(G/B, \L({\lambda})) \rightarrow H^0(G/B, \L(V_1)) \\
&\rightarrow H^1(G/B, \L({\lambda_1})) \rightarrow
H^{p,1}(G/B, \L({\lambda})) \rightarrow H^1(G/B,\L(V_1))\\
&\rightarrow H^2(G/B, \L ({\lambda_1}))\rightarrow H^{p,2}(G/B,
\L({\lambda})) \rightarrow H^2(G/B, \L(V_1))\rightarrow \cdots
\end{align*}

By our arrangement of weights $\lambda_2$ will be the lowest weight
among all the weights of $V_1$. Let $V_2 = V_1/k_{\lambda_2}$. As before
we will have a short exact sequence $0\rightarrow k_{\lambda_1}
\rightarrow V_1 \rightarrow V_2 \rightarrow0$. Which will induce the
long exact sequences of cohomology modules. Proceeding like this we will
have several exact sequences, but by Lemma 10, we essentially need to
consider only three (in case of $G_2$ four) long exact sequences which
involve the weights $\mu_1, ~ \mu_2, ~0$. The three exact sequence are
as follows:
\begin{align*}
0 &\rightarrow H^0(G/B,\L({\mu_1})) \rightarrow H^0(G/B, \L(V_s))
\rightarrow H^0(G/B,\L(V_{s+1}))\\
&\rightarrow H^1(G/B, \L({\mu_1})) \rightarrow H^1(G/B, \L(V_s))
\rightarrow H^1(G/B,\L(V_{s+1}))\\
&\rightarrow H^2(G/B, \L (\mu_1))
\rightarrow H^2(G/B, \L(V_s)) \rightarrow \cdots\\
0 &\rightarrow H^0(G/B,\L({\mu_2})) \rightarrow H^0(G/B, \L(V_{s+1}))
\rightarrow H^0(G/B, \L(V_{s+2}))\\
&\rightarrow H^1(G/B, \L({\mu_2}))
\rightarrow H^1(G/B, \L(V_{s+1})) \rightarrow H^1(G/B,\L(V_{s+2}))\\
&\rightarrow H^2(G/B, \L (\mu_2))\rightarrow \cdots\\
0 &\rightarrow H^0(G/B,\L(0)) \rightarrow H^0(G/B, \L(V_t))
\rightarrow H^0(G/B, \L(V_{t+1}))\\
&\rightarrow H^1(G/B, \L(0)) \rightarrow H^1(G/B, \L(V_t))
\rightarrow H^1(G/B,\L(V_{t+1}))\\
&\rightarrow H^2(G/B, \L (0))\rightarrow H^2(G/B, \L(V_t))
\rightarrow  \cdots
\end{align*}

If we use $(\star)$ and the Borel--Weil--Bott theorem, we could easily
see that $H^1(G/B, \L(V_s))$ $\neq 0$. Again by Lemma 10 and the fact that
weights appeared before $\mu_1$ are all in $X(T)-W\cdot X(T)^+$, we can
easily see that $H^{p,1}(G/B, \L({\lambda}))\neq0$. \hfill $\BBox\,$
\end{proof}

The following corollary is an immediate consequence of Theorem~12.

\setcounter{coro}{12}
\begin{coro}$\left.\right.$\vspace{.5pc}

\noindent Let $G$ be a simple simply connected algebraic group with {\rm
rank}$(G)\geq 2$. Then $G/B$ does not satisfy Bott vanishing property.\vspace{-.2pc}
\end{coro}

\setcounter{rem}{13}
\begin{rem}
{\rm Our motivation to prove Theorem~11 was Corollary~13. This is the
reason why we wanted strict dominant weights. If we just want the result
for dominant weights there are families of weights for which higher
cohomology survives. For example, let $G$ be of type $A_2$. Let $\lambda
= n\omega_1+ m\omega_2$. If $n=0$ and $m\geq 2$ or $n\geq 2$ and $m=0$
then $H^{2,1}(G/B, \L({\lambda}))\not = 0$.}\vspace{-.2pc}
\end{rem}

\begin{rem}
{\rm Let $\lambda = \sum n_k \omega_k$ and $d=\dim(G/B)$. If $n_k \geq 1,
\forall k$ then $L_{\lambda}$ will be ample. By Theorem 1, for $i\geq 1$
to get $H^{d-1,i}(G/B, \L_{\lambda})\not =0$ we cannot assume all $n_k$
are `bigger'. One natural candidate with all its coefficients positive
is $\rho = ( 1, \ldots , 1).$ But we note that if $G$ is not of type
$A_2$ all cohomology groups vanish for $\lambda=\rho$.}\vspace{-.2pc}
\end{rem}

\setcounter{theorr}{15}
\begin{therr}[\!]
Let $G$ be a simple{\rm ,} simply connected algebraic group with {\rm
rank}$(G) \geq 2$. Then\vspace{-.3pc}
\begin{enumerate}
\renewcommand\labelenumi{\rm (\roman{enumi})}
\leftskip .15pc
\item $G/B$ is not a toric variety.
\item $G/B$ cannot be degenerated to a smooth toric variety in such a
manner that ample cone degenerates to ample cone.\vspace{-1pc}
\end{enumerate}
\end{therr}

\begin{proof}$\left.\right.$\vspace{-.3pc}
\begin{enumerate}
\renewcommand\labelenumi{\rm (\roman{enumi})}
\leftskip .15pc
\item Since $G/B$ is smooth, if it is a toric variety it has to satisfy
Bott vanishing property (cf. \cite{btlm}), which contradicts Corollary~13.

\item Suppose it degenerates to a smooth toric variety such that the
ample cone degenerates to the ample cone, then by \cite{btlm} and
semicontinuity theorem, $G/B$ will satisfy Bott vanishing property,
which contradicts Corollary~13.\hfill $\BBox\ \,$\vspace{-1.5pc}
\end{enumerate}
\end{proof}

\setcounter{rem}{16}
\begin{rem}$\left.\right.$
{\rm \begin{enumerate}
\renewcommand{\labelenumi}{\arabic{enumi}.}
\leftskip -.2pc
\item The fact that $G/B$ is not a toric variety could possibly be
checked by pure topological criteria, for example by computing the
integral cohomology and observing that for smooth projective toric
variety the integral cohomology ring is generated by the second
cohomology. Since the integral cohomology of $SL(n)/B$ is generated as a
ring by degree 2 elements, this criterion will not help to conclude
$SL(n)/B$ is not a toric\break variety.

\item One could compare the rank of the connected component of
automorphism group of $G/B(= G/Z(G)$ where $Z(G)$ is the center of $G$)
with the dimension of $G/B$.

\item The above methods can be used to prove Theorem~16(i) but not
Corollary~13 or Theorem~16(ii).
\end{enumerate}}
\end{rem}

\section*{Appendix A}

The following is the set of all positive roots of exceptional groups. We
are presenting them in terms of both simple roots and fundamental
weights. It will be useful to order them and to apply Borel--Weil--Bott
theorem. The left-hand side is the coefficient of positive roots with
respect to the simple roots, and the right-hand side is with respect to
the fundamental weights.

\renewcommand\thesubsection{A.\arabic{subsection}}

\setcounter{subsection}{0}
\subsection{\it Positive roots of $G_2$ and $F_4$}

{\footnotesize
$$      \begin{tabular}{@{}lc@{\ \ \ \ \ \ \ \ \ \ }rr r lc@{\ \ \ \ \ \ \ \ \ \ }rrrr@{}}
& & ${\boldmath G_2}$ & & & & & ${\boldmath F_4}$ & & & \\
& &  & & & & &  & & & \\
(1 0)& \myfavouritearrow & ($\phantom{-}$2 & $-$3)& & (1 0 0 0) &
\myfavouritearrow & ($\phantom{-}$2 & $-$1& 0 & 0) \\
(0 1)& \myfavouritearrow & ($-$1 & 2)& & (0 1 0 0)&
\myfavouritearrow & ($-$1 & 2& $-$2 & 0) \\
(1 1)& \myfavouritearrow & ($\phantom{-}$1 & $-$1)& & (0 0 1 0)&
\myfavouritearrow & ($\phantom{-}$0 & $-$1& 2 & $-$1) \\
(1 2)& \myfavouritearrow & ($\phantom{-}$0 & 1)& & (0 0 0 1)&
\myfavouritearrow & ($\phantom{-}$0 &  0&$-$1 & 2) \\
(1 3)& \myfavouritearrow & ($-$1 & 3)& & (1 1 0 0)&
\myfavouritearrow & ($\phantom{-}$1 &  1& $-$2  & 0) \\
(2 3)& \myfavouritearrow & ($\phantom{-}$1 & 0)& & (0 1 1 0)&
\myfavouritearrow & ($-$1 & 1& 0 &$-$1) \\
& & & & & (0 0 1 1)& \myfavouritearrow & ($\phantom{-}$0 & $-$1 & 1 & 1)\\
& & & & & (1 1 1 0)& \myfavouritearrow & ($\phantom{-}$1 & 0 & 0 & $-$1)\\
& & & & & (0 1 1 1)& \myfavouritearrow & ($-$1 & 1 & $-$1 & 1)\\
& & & & & (1 1 l 1)& \myfavouritearrow & ($\phantom{-}$1 & 0    & $-$1 & 1)\\
& & & & & (0 1 2 0)& \myfavouritearrow & ($-$1 & 0 & 2 & $-$2)\\
& & & & & (1 1 2 0)& \myfavouritearrow & ($\phantom{-}$1 & $-$1 & 2  & $-$2)\\
& & & & & (0 1 2 1)& \myfavouritearrow & ($-$1 & 0 & 1 & 0)\\
& & & & & (1 2 2 0)& \myfavouritearrow & ($\phantom{-}$0 & 1 & 0   & $-$2)\\
& & & & & (1 1 2 1)& \myfavouritearrow & ($\phantom{-}$1 & $-$1 & 1 & 0)\\
& & & & & (0 1 2 2)& \myfavouritearrow & ($-$1 & 0 & 0 & 2)\\
& & & & & (1 1 2 2)& \myfavouritearrow & ($\phantom{-}$1 & $-$1 & 0   & 2)\\
& & & & & (1 2 2 1)& \myfavouritearrow & ($\phantom{-}$0 & 1 & $-$1 & 0)\\
& & & & & (1 2 2 2)& \myfavouritearrow & ($\phantom{-}$0 & 1 & $-$2 & 2)\\
& & & & & (1 2 3 1)& \myfavouritearrow & ($\phantom{-}$0 & 0 & 1 & $-$1)\\
& & & & & (1 2 3 2)& \myfavouritearrow & ($\phantom{-}$0 & 0 & 0   & 1)\\
& & & & & (1 2 4 2)& \myfavouritearrow & ($\phantom{-}$0 & $-$1 & 2 & 0)\\
& & & & & (1 3 4 2)& \myfavouritearrow & ($-$1 & 1 & 0 & 0)\\
& & & & & (2 3 4 2)& \myfavouritearrow & ($\phantom{-}$1 & 0 & 0   & 0)\\
\end{tabular} $$}

\subsection{\it Positive roots of $E_6$}

{\footnotesize
$$\begin{tabular}{lc@{\ \ \ \ \ \ \ \ \ \ }rrrrrr}
(1  0  0  0  0  0) & \myfavouritearrow &  ($\phantom{-}$2 & 0& $-$1&  0 & 0 & 0)\\
(0  1  0  0  0  0) & \myfavouritearrow &  ($\phantom{-}$0&  2&  0& $-$1&  0 & 0)\\
(0  0  1  0  0  0) & \myfavouritearrow &  ($-$1 & 0 & 2 &$-$1 & 0  &0)\\
(0  0  0  1  0  0) & \myfavouritearrow &  ($\phantom{-}$0& $-$1& $-$1&  2& $-$1&  0)\\
(0  0  0  0  1  0) & \myfavouritearrow &  ($\phantom{-}$0 & 0 & 0 &$-$1 & 2 &$-$1)\\
(0  0  0  0  0  1) & \myfavouritearrow &  ($\phantom{-}$0&  0&  0 & 0& $-$1&  2)\\
(1  0  1  0  0  0) & \myfavouritearrow &  ($\phantom{-}$1 & 0 & 1 &$-$1 & 0 & 0)\\
(0  1  0  1  0  0) & \myfavouritearrow &  ($\phantom{-}$0&  1 &$-$1 & 1 &$-$1 & $-$1)\\
(0  0  1  1  0  0) & \myfavouritearrow &  ($-$1& $-$1 & 1 & 1 &$-$1 & 0)\\
(0  0  0  1  1  0) & \myfavouritearrow &  ($\phantom{-}$0& $-$1 &$-$1 & 1 & 1 &$-$1)\\
(0  0  0  0  1  1) & \myfavouritearrow &  ($\phantom{-}$0&  0 & 0 &$-$1 & 1 & 1)\\
(1  0  1  1  0  0) & \myfavouritearrow &  ($\phantom{-}$1& $-$1 & 0 & 1 &$-$1 & 0)\\
(0  1  1  1  0  0) & \myfavouritearrow &  ($-$1&  1 & 1 & 0 &$-$1 & 0)\\
(0  1  0  1  1  0) & \myfavouritearrow &  ($\phantom{-}$0&  1 &$-$1 & 0 & 1 &$-$1)\\
(0  0  1  1  1  0) & \myfavouritearrow &  ($-$1& $-$1 & 1 & 0 & 1 &$-$1)\\
(0  0  0  1  1  1) & \myfavouritearrow & ($\phantom{-}$0& $-$1 &$-$1 & 1 & 0 & 1)\\
(1  1  1  1  0  0) & \myfavouritearrow & ($\phantom{-}$1&  1 & 0 & 0 &$-$1 & 0)\\
(1  0  1  1  1  0) & \myfavouritearrow & ($\phantom{-}$1& $-$1 & 0 & 0 & 1 &$-$1)\\
(0  1  1  1  1  0) & \myfavouritearrow & ($-$1&  1 & 1 &$-$1 & 1 &$-$1)\\
(0  1  0  1  1  1) & \myfavouritearrow & ($\phantom{-}$0&  1 &$-$1 & 0 & 0 & 1)\\
(0  0  1  1  1  1) & \myfavouritearrow & ($-$1& $-$1 & 1 & 0 & 0 & 1)\\
(1  1  1  1  1  0) & \myfavouritearrow & ($\phantom{-}$1&  1 & 0 &$-$1 & 1 &$-$1)\\
(1  0  1  1  1  1) & \myfavouritearrow & ($\phantom{-}$1& $-$1 & 0 & 0 & 0 & 1)\\
(0  1  1  1  1  1) & \myfavouritearrow & ($-$1&  1 & 1 &$-$1 & 0 & 1)\\
(0  1  1  2  1  0) & \myfavouritearrow & ($-$1&  0 & 0 & 1 & 0 &$-$1)\\
(1  1  1  1  1  1) & \myfavouritearrow & ($\phantom{-}$1&  1 & 0 &$-$1 & 0 & 1)\\
(0  1  1  2  1  1) & \myfavouritearrow & ($-$1&  0 & 0 & 1 &$-$1 & 1)\\
(1  1  1  2  1  0) & \myfavouritearrow & ($\phantom{-}$1&  0 &$-$1 & 1 & 0 &$-$1)\\
(1  1  1  2  1  1) & \myfavouritearrow & ($\phantom{-}$1&  0 &$-$1 & 1 &$-$1 & 1)\\
(0  1  1  2  2  1) & \myfavouritearrow & ($-$1&  0 & 0 & 0 & 1 & 0)\\
(1  1  2  2  1  0) & \myfavouritearrow & ($\phantom{-}$0&  0 & 1 & 0 & 0 &$-$1)\\
(1  1  1  2  2  1) & \myfavouritearrow & ($\phantom{-}$1&  0 &$-$1 & 0 & 1 & 0)\\
(1  1  2  2  1  1) & \myfavouritearrow & ($\phantom{-}$0&  0 & 1 & 0 &$-$1 & 1)\\
(1  1  2  2  2  1) & \myfavouritearrow & ($\phantom{-}$0&  0 & 1 &$-$1 & 1 & 0)\\
(1  1  2  3  2  1) & \myfavouritearrow & ($\phantom{-}$0& $-$1 & 0 & 1 & 0 & 0)\\
(1  2  2  3  2  1) & \myfavouritearrow & ($\phantom{-}$0&  1 & 0 & 0 & 0 & 0)%\\
\end{tabular}$$}
\subsection{\it Positive roots of $E_7$}
{\footnotesize
$$\begin{tabular}{lc@{\ \ \ \ \ \ \ \ \ \ }rrrrrrr}
(1  0  0  0  0  0  0) & \myfavouritearrow & ($\phantom{-}$2 & 0 &$-$1 & 0 & 0 & 0 & 0)\\
(0  1  0  0  0  0  0) & \myfavouritearrow & ($\phantom{-}$0 & 2 & 0 &$-$1 & 0 & 0 & 0)\\
(0  0  1  0  0  0  0) & \myfavouritearrow & ($-$1 & 0 & 2 &$-$1 & 0 & 0 & 0)\\
(0  0  0  1  0  0  0) & \myfavouritearrow & ($\phantom{-}$0 &$-$1 &$-$1 & 2 &$-$1 & 0 & 0)\\
(0  0  0  0  1  0  0) & \myfavouritearrow & ($\phantom{-}$0 & 0 & 0 &$-$1 & 2 &$-$1 & 0)\\
(0  0  0  0  0  1  0) & \myfavouritearrow & ($\phantom{-}$0 & 0 & 0 & 0 &$-$1 & 2 &$-$1)\\
(0  0  0  0  0  0  1) & \myfavouritearrow & ($\phantom{-}$0 & 0 & 0 & 0 & 0 &$-$1 & 2)\\
(1  0  1  0  0  0  0) & \myfavouritearrow & ($\phantom{-}$1 & 0 & 1 &$-$1 & 0 & 0 & 0)\\
(0  1  0  1  0  0  0) & \myfavouritearrow & ($\phantom{-}$0 & 1 &$-$1 & 1 &$-$1 & 0 & 0)\\
(0  0  1  1  0  0  0) & \myfavouritearrow & ($-$1 &$-$1 & 1 & 1 &$-$1 & 0 & 0)\\
(0  0  0  1  1  0  0) & \myfavouritearrow & ($\phantom{-}$0 &$-$1 &$-$1 & 1 & 1 &$-$1 & 0)\\
(0  0  0  0  1  1  0) & \myfavouritearrow & ($\phantom{-}$0 & 0 & 0 &$-$1 & 1 & 1 &$-$1)\\
(0  0  0  0  0  1  1) & \myfavouritearrow & ($\phantom{-}$0 & 0 & 0 & 0 &$-$1 & 1 & 1)\\
(1  0  1  1  0  0  0) & \myfavouritearrow & ($\phantom{-}$1 &$-$1 & 0 & 1 &$-$1 & 0 & 0)%\\
\end{tabular}$$

$$\begin{tabular}{lc@{\ \ \ \ \ \ \ \ \ \ }rrrrrrr}
(0  1  1  1  0  0  0) & \myfavouritearrow & ($-$1 & 1 & 1 & 0 &$-$1 & 0 & 0)\\
(0  1  0  1  1  0  0) & \myfavouritearrow & ($\phantom{-}$0 & 1 &$-$1 & 0 & 1 &$-$1 & 0)\\
(0  0  1  1  1  0  0) & \myfavouritearrow & ($-$1 &$-$1 & 1 & 0 & 1 &$-$1 & 0)\\
(0  0  0  1  1  1  0) & \myfavouritearrow & ($\phantom{-}$0 &$-$1 &$-$1 & 1 & 0 & 1 &$-$1)\\
(0  0  0  0  1  1  1) & \myfavouritearrow & ($\phantom{-}$0 & 0 & 0 &$-$1 & 1 & 0 & 1)\\
(1  1  1  1  0  0  0) & \myfavouritearrow & ($\phantom{-}$1 & 1 & 0 & 0 &$-$1 & 0 & 0)\\
(1  0  1  1  1  0  0) & \myfavouritearrow & ($\phantom{-}$1 &$-$1 & 0 & 0 & 1 &$-$1 & 0)\\
(0  1  1  1  1  0  0) & \myfavouritearrow & ($-$1 & 1 & 1 &$-$1 & 1 &$-$1 & 0)\\
(0  1  0  1  1  1  0) & \myfavouritearrow & ($\phantom{-}$0 & 1 &$-$1 & 0 & 0 & 1 &$-$1)\\
(0  0  1  1  1  1  0) & \myfavouritearrow & ($-$1 &$-$1 & 1 & 0 & 0 & 1 &$-$1)\\
(0  0  0  1  1  1  1) & \myfavouritearrow & ($\phantom{-}$0 &$-$1 &$-$1 & 1 & 0 & 0 & 1)\\
(1  1  1  1  1  0  0) & \myfavouritearrow & ($\phantom{-}$1 & 1 & 0 &$-$1 & 1 &$-$1 & 0)\\
(1  0  1  1  1  1  0) & \myfavouritearrow & ($\phantom{-}$1 &$-$1 & 0 & 0 & 0 & 1 &$-$1)\\
(0  1  1  1  1  1  0) & \myfavouritearrow & ($-$1 & 1 & 1 &$-$1 & 0 & 1 &$-$1)\\
(0  1  0  1  1  1  1) & \myfavouritearrow & ($\phantom{-}$0 & 1 &$-$1 & 0 & 0 & 0 & 1)\\
(0  0  1  1  1  1  1) & \myfavouritearrow & ($-$1 &$-$1 & 1 & 0 & 0 & 0 & 1)\\
(0  1  1  2  1  0  0) & \myfavouritearrow & ($-$1 & 0 & 0 & 1 & 0 &$-$1 & 0)\\
(1  1  1  1  1  1  0) & \myfavouritearrow & ($\phantom{-}$1 & 1 & 0 &$-$1 & 0 & 1 &$-$1)\\
(1  0  1  1  1  1  1) & \myfavouritearrow & ($\phantom{-}$1 &$-$1 & 0 & 0 & 0 & 0 & 1)\\
(0  1  1  1  1  1  1) & \myfavouritearrow & ($-$1 & 1 & 1 &$-$1 & 0 & 0 & 1)\\
(0  1  1  2  1  1  0) & \myfavouritearrow & ($-$1 & 0 & 0 & 1 &$-$1 & 1 &$-$1)\\
(1  1  1  2  1  0  0) & \myfavouritearrow & ($\phantom{-}$1 & 0 &$-$1 & 1 & 0 &$-$1 & 0)\\
(1  1  1  1  1  1  1) & \myfavouritearrow & ($\phantom{-}$1 & 1 & 0 &$-$1 & 0 & 0 & 1)\\
(1  1  1  2  1  1  0) & \myfavouritearrow & ($\phantom{-}$1 & 0 &$-$1 & 1 &$-$1 & 1 &$-$1)\\
(0  1  1  2  1  1  1) & \myfavouritearrow & ($-$1 & 0 & 0 & 1 &$-$1 & 0 & 1)\\
(0  1  1  2  2  1  0) & \myfavouritearrow & ($-$1 & 0 & 0 & 0 & 1 & 0 &$-$1)\\
(1  1  2  2  1  0  0) & \myfavouritearrow & ($\phantom{-}$0 & 0 & 1 & 0 & 0 &$-$1 & 0)\\
(1  1  1  2  1  1  1) & \myfavouritearrow & ($\phantom{-}$1 & 0 &$-$1 & 1 &$-$1 & 0 & 1)\\
(1  1  1  2  2  1  0) & \myfavouritearrow & ($\phantom{-}$1 & 0 &$-$1 & 0 & 1 & 0 &$-$1)\\
(1  1  2  2  1  1  0) & \myfavouritearrow & ($\phantom{-}$0 & 0 & 1 & 0 &$-$1 & 1 &$-$1)\\
(0  1  1  2  2  1  1) & \myfavouritearrow & ($-$1 & 0 & 0 & 0 & 1 &$-$1 & 1)\\
(1  1  2  2  1  1  1) & \myfavouritearrow & ($\phantom{-}$0 & 0 & 1 & 0 &$-$1 & 0 & 1)\\
(1  1  2  2  2  1  0) & \myfavouritearrow & ($\phantom{-}$0 & 0 & 1 &$-$1 & 1 & 0 &$-$1)\\
(1  1  1  2  2  1  1) & \myfavouritearrow & ($\phantom{-}$1 & 0 &$-$1 & 0 & 1 &$-$1 & 1)\\
(0  1  1  2  2  2  1) & \myfavouritearrow & ($-$1 & 0 & 0 & 0 & 0 & 1 & 0)\\
(1  1  2  2  2  1  1) & \myfavouritearrow & ($\phantom{-}$0 & 0 & 1 &$-$1 & 1 &$-$1 & 1)\\
(1  1  1  2  2  2  1) & \myfavouritearrow & ($\phantom{-}$1 & 0 &$-$1 & 0 & 0 & 1 & 0)\\
(1  1  2  3  2  1  0) & \myfavouritearrow & ($\phantom{-}$0 &$-$1 & 0 & 1 & 0 & 0 &$-$1)\\
(1  1  2  2  2  2  1) & \myfavouritearrow & ($\phantom{-}$0 & 0 & 1 &$-$1 & 0 & 1 & 0)\\
(1  2  2  3  2  1  0) & \myfavouritearrow & ($\phantom{-}$0 & 1 & 0 & 0 & 0 & 0 &$-$1)\\
(1  1  2  3  2  1  1) & \myfavouritearrow & ($\phantom{-}$0 &$-$1 & 0 & 1 & 0 &$-$1 & 1)\\
(1  1  2  3  2  2  1) & \myfavouritearrow & ($\phantom{-}$0 &$-$1 & 0 & 1 &$-$1 & 1 & 0)\\
(1  2  2  3  2  1  1) & \myfavouritearrow & ($\phantom{-}$0 & 1 & 0 & 0 & 0 &$-$1 & 1)\\
(1  1  2  3  3  2  1) & \myfavouritearrow & ($\phantom{-}$0 &$-$1 & 0 & 0 & 1 & 0 & 0)\\
(1  2  2  3  2  2  1) & \myfavouritearrow & ($\phantom{-}$0 & 1 & 0 & 0 &$-$1 & 1 & 0)\\
(1  2  2  3  3  2  1) & \myfavouritearrow & ($\phantom{-}$0 & 1 & 0 &$-$1 & 1 & 0 & 0)\\
(1  2  2  4  3  2  1) & \myfavouritearrow & ($\phantom{-}$0 & 0 &$-$1 & 1 & 0 & 0 & 0)\\
(1  2  3  4  3  2  1) & \myfavouritearrow & ($-$1 & 0 & 1 & 0 & 0 & 0 & 0)\\
(2  2  3  4  3  2  1) & \myfavouritearrow & ($\phantom{-}$1 & 0 & 0 & 0 & 0 & 0 & 0)%\\
\end{tabular}$$}

\subsection{\it Positive roots of $E_8$}

{\footnotesize
$$\begin{tabular}{lc@{\ \ \ \ \ \ \ \ \ \ }rrrrrrrr}
(1  0  0  0  0  0  0  0) & \myfavouritearrow & ($\phantom{-}$2&  0& $-$1&  0&  0&  0&  0&  0)\\
(0  1  0  0  0  0  0  0) & \myfavouritearrow & ($\phantom{-}$0&  2&  0& $-$1&  0&  0&  0&  0)\\
(0  0  1  0  0  0  0  0) & \myfavouritearrow & ($-$1&  0&  2& $-$1&  0&  0&  0&  0)\\
(0  0  0  1  0  0  0  0) & \myfavouritearrow & ($\phantom{-}$0& $-$1& $-$1&  2& $-$1&  0&  0&  0)%\\
\end{tabular}$$

$$\begin{tabular}{lc@{\ \ \ \ \ \ \ \ \ \ }rrrrrrrr}
(0  0  0  0  1  0  0  0) & \myfavouritearrow & ($\phantom{-}$0&  0&  0& $-$1&  2& $-$1&  0&  0)\\
(0  0  0  0  0  1  0  0) & \myfavouritearrow & ($\phantom{-}$0&  0&  0&  0& $-$1&  2& $-$1&  0)\\
(0  0  0  0  0  0  1  0) & \myfavouritearrow & ($\phantom{-}$0&  0&  0&  0&  0& $-$1&  2& $-$1)\\
(0  0  0  0  0  0  0  1) & \myfavouritearrow & ($\phantom{-}$0&  0&  0&  0&  0&  0& $-$1&  2)\\
(1  0  1  0  0  0  0  0) & \myfavouritearrow & ($\phantom{-}$1&  0&  1& $-$1&  0&  0&  0&  0)\\
(0  1  0  1  0  0  0  0) & \myfavouritearrow & ($\phantom{-}$0&  1& $-$1&  1& $-$1&  0&  0&  0)\\
(0  0  1  1  0  0  0  0) & \myfavouritearrow & ($-$1& $-$1&  1&  1& $-$1&  0&  0&  0)\\
(0  0  0  1  1  0  0  0) & \myfavouritearrow & ($\phantom{-}$0& $-$1& $-$1&  1&  1& $-$1&  0&  0)\\
(0  0  0  0  1  1  0  0) & \myfavouritearrow & ($\phantom{-}$0&  0&  0& $-$1&  1&  1& $-$1&  0)\\
(0  0  0  0  0  1  1  0) & \myfavouritearrow & ($\phantom{-}$0&  0&  0&  0& $-$1&  1&  1& $-$1)\\
(0  0  0  0  0  0  1  1) & \myfavouritearrow & ($\phantom{-}$0&  0&  0&  0&  0& $-$1&  1&  1)\\
(1  0  1  1  0  0  0  0) & \myfavouritearrow & ($\phantom{-}$1& $-$1&  0&  1& $-$1&  0&  0&  0)\\
(0  1  1  1  0  0  0  0) & \myfavouritearrow & ($-$1&  1&  1&  0& $-$1&  0&  0&  0)\\
(0  1  0  1  1  0  0  0) & \myfavouritearrow & ($\phantom{-}$0&  1& $-$1&  0&  1& $-$1&  0&  0)\\
(0  0  1  1  1  0  0  0) & \myfavouritearrow & ($-$1& $-$1&  1&  0&  1& $-$1&  0&  0)\\
(0  0  0  1  1  1  0  0) & \myfavouritearrow & ($\phantom{-}$0& $-$1& $-$1&  1&  0&  1& $-$1&  0)\\
(0  0  0  0  1  1  1  0) & \myfavouritearrow & ($\phantom{-}$0&  0&  0& $-$1&  1&  0&  1& $-$1)\\
(0  0  0  0  0  1  1  1) & \myfavouritearrow & ($\phantom{-}$0&  0&  0&  0& $-$1&  1&  0&  1)\\
(1  1  1  1  0  0  0  0) & \myfavouritearrow & ($\phantom{-}$1&  1&  0&  0& $-$1&  0&  0&  0)\\
(1  0  1  1  1  0  0  0) & \myfavouritearrow & ($\phantom{-}$1& $-$1&  0&  0&  1& $-$1&  0&  0)\\
(0  1  1  1  1  0  0  0) & \myfavouritearrow & ($-$1&  1&  1& $-$1&  1& $-$1&  0&  0)\\
(0  1  0  1  1  1  0  0) & \myfavouritearrow & ($\phantom{-}$0&  1& $-$1&  0&  0&  1& $-$1&  0)\\
(0  0  1  1  1  1  0  0) & \myfavouritearrow & ($-$1& $-$1&  1&  0&  0&  1& $-$1&  0)\\
(0  0  0  1  1  1  1  0) & \myfavouritearrow & ($\phantom{-}$0& $-$1& $-$1&  1&  0&  0&  1& $-$1)\\
(0  0  0  0  1  1  1  1) & \myfavouritearrow & ($\phantom{-}$0&  0&  0& $-$1&  1&  0&  0&  1)\\
(1  1  1  1  1  0  0  0) & \myfavouritearrow & ($\phantom{-}$1&  1&  0& $-$1&  1& $-$1&  0&  0)\\

(1  0  1  1  1  1  0  0) & \myfavouritearrow & ($\phantom{-}$1& $-$1&  0&  0&  0&  1& $-$1&  0)\\
(0  1  1  1  1  1  0  0) & \myfavouritearrow & ($-$1&  1&  1& $-$1&  0&  1& $-$1&  0)\\
(0  1  0  1  1  1  1  0) & \myfavouritearrow & ($\phantom{-}$0&  1& $-$1&  0&  0&  0&  1& $-$1)\\
(0  0  1  1  1  1  1  0) & \myfavouritearrow & ($-$1& $-$1&  1&  0&  0&  0&  1& $-$1)\\

(0  0  0  1  1  1  1  1) & \myfavouritearrow & ($\phantom{-}$0& $-$1& $-$1&  1&  0&  0&  0&  1)\\
(0  1  1  2  1  0  0  0) & \myfavouritearrow & ($-$1&  0&  0&  1&  0& $-$1&  0&  0)\\
(1  1  1  1  1  1  0  0) & \myfavouritearrow & ($\phantom{-}$1&  1&  0& $-$1&  0&  1& $-$1&  0)\\
(1  0  1  1  1  1  1  0) & \myfavouritearrow & ($\phantom{-}$1& $-$1&  0&  0&  0&  0&  1& $-$1)\\
(0  1  1  1  1  1  1  0) & \myfavouritearrow & ($-$1&  1&  1& $-$1&  0&  0&  1& $-$1)\\

(0  1  0  1  1  1  1  1) & \myfavouritearrow & ($\phantom{-}$0&  1& $-$1&  0&  0&  0&  0&  1)\\
(0  0  1  1  1  1  1  1) & \myfavouritearrow & ($-$1& $-$1&  1&  0&  0&  0&  0&  1)\\
(0  1  1  2  1  1  0  0) & \myfavouritearrow & ($-$1&  0&  0&  1& $-$1&  1& $-$1&  0)\\
(1  1  1  2  1  0  0  0) & \myfavouritearrow & ($\phantom{-}$1&  0& $-$1&  1&  0& $-$1&  0&  0)\\
(1  1  1  1  1  1  1  0) & \myfavouritearrow & ($\phantom{-}$1&  1&  0& $-$1&  0&  0&  1& $-$1)\\
(1  0  1  1  1  1  1  1) & \myfavouritearrow & ($\phantom{-}$1& $-$1&  0&  0&  0&  0&  0&  1)\\
(0  1  1  1  1  1  1  1) & \myfavouritearrow & ($-$1&  1&  1& $-$1&  0&  0&  0&  1)\\

(1  1  1  2  1  1  0  0) & \myfavouritearrow & ($\phantom{-}$1&  0& $-$1&  1& $-$1&  1& $-$1&  0)\\
(0  1  1  2  1  1  1  0) & \myfavouritearrow & ($-$1&  0&  0&  1& $-$1&  0&  1& $-$1)\\
(0  1  1  2  2  1  0  0) & \myfavouritearrow & ($-$1&  0&  0&  0&  1&  0& $-$1&  0)\\
(1  1  2  2  1  0  0  0) & \myfavouritearrow & ($\phantom{-}$0&  0&  1&  0&  0& $-$1&  0&  0)\\
(1  1  1  1  1  1  1  1) & \myfavouritearrow & ($\phantom{-}$1&  1&  0& $-$1&  0&  0&  0&  1)\\
(1  1  1  2  1  1  1  0) & \myfavouritearrow & ($\phantom{-}$1&  0& $-$1&  1& $-$1&  0&  1& $-$1)\\
(1  1  1  2  2  1  0  0) & \myfavouritearrow & ($\phantom{-}$1&  0& $-$1&  0&  1&  0& $-$1&  0)\\
(1  1  2  2  1  1  0  0) & \myfavouritearrow & ($\phantom{-}$0&  0&  1&  0& $-$1&  1& $-$1&  0)\\
(0  1  1  2  1  1  1  1) & \myfavouritearrow & ($-$1&  0&  0&  1& $-$1&  0&  0&  1)\\
(0  1  1  2  2  1  1  0) & \myfavouritearrow & ($-$1&  0&  0&  0&  1& $-$1&  1& $-$1)\\
(1  1  1  2  1  1  1  1) & \myfavouritearrow & ($\phantom{-}$1&  0& $-$1&  1& $-$1&  0&  0&  1)\\
(0  1  1  2  2  1  1  1) & \myfavouritearrow & ($-$1&  0&  0&  0&  1& $-$1&  0&  1)\\
(1  1  2  2  1  1  1  0) & \myfavouritearrow & ($\phantom{-}$0&  0&  1&  0& $-$1&  0&  1& $-$1)\\
(1  1  2  2  2  1  0  0) & \myfavouritearrow & ($\phantom{-}$0&  0&  1& $-$1&  1&  0& $-$1&  0)\\
(1  1  1  2  2  1  1  0) & \myfavouritearrow & ($\phantom{-}$1&  0& $-$1&  0&  1& $-$1&  1& $-$1)%\\
\end{tabular}$$

$$\begin{tabular}{lc@{\ \ \ \ \ \ \ \ \ \ }rrrrrrrr}
(0  1  1  2  2  2  1  0) & \myfavouritearrow & ($-$1&  0&  0&  0&  0&  1&  0& $-$1)\\[.1pc]
(0  1  1  2  2  2  1  1) & \myfavouritearrow & ($-$1&  0&  0&  0&  0&  1& $-$1&  1)\\[.1pc]
(1  1  2  2  1  1  1  1) & \myfavouritearrow & ($\phantom{-}$0&  0&  1&  0& $-$1&  0&  0&  1)\\[.1pc]
(1  1  1  2  2  1  1  1) & \myfavouritearrow & ($\phantom{-}$1&  0& $-$1&  0&  1& $-$1&  0&  1)\\[.1pc]
(1  1  2  2  2  1  1  0) & \myfavouritearrow & ($\phantom{-}$0&  0&  1& $-$1&  1& $-$1&  1& $-$1)\\[.1pc]
(1  1  1  2  2  2  1  0) & \myfavouritearrow & ($\phantom{-}$1&  0& $-$1&  0&  0&  1&  0& $-$1)\\[.1pc]
(1  1  2  3  2  1  0  0) & \myfavouritearrow & ($\phantom{-}$0& $-$1&  0&  1&  0&  0& $-$1&  0)\\[.1pc]
(1  1  2  2  2  1  1  1) & \myfavouritearrow & ($\phantom{-}$0&  0&  1& $-$1&  1& $-$1&  0&  1)\\[.1pc]
(1  1  1  2  2  2  1  1) & \myfavouritearrow & ($\phantom{-}$1&  0& $-$1&  0&  0&  1& $-$1&  1)\\[.1pc]
(0  1  1  2  2  2  2  1) & \myfavouritearrow & ($-$1&  0&  0&  0&  0&  0&  1&  0)\\[.1pc]
(1  1  2  2  2  2  1  0) & \myfavouritearrow & ($\phantom{-}$0&  0&  1& $-$1&  0&  1&  0& $-$1)\\[.1pc]
(1  2  2  3  2  1  0  0) & \myfavouritearrow & ($\phantom{-}$0&  1&  0&  0&  0&  0& $-$1&  0)\\[.1pc]
(1  1  2  3  2  1  1  0) & \myfavouritearrow & ($\phantom{-}$0& $-$1&  0&  1&  0& $-$1&  1& $-$1)\\[.1pc]
(1  1  2  3  2  1  1  1) & \myfavouritearrow & ($\phantom{-}$0& $-$1&  0&  1&  0& $-$1&  0&  1)\\[.1pc]
(1  1  2  2  2  2  1  1) & \myfavouritearrow & ($\phantom{-}$0&  0&  1& $-$1&  0&  1& $-$1&  1)\\[.1pc]
(1  1  1  2  2  2  2  1) & \myfavouritearrow & ($\phantom{-}$1&  0& $-$1&  0&  0&  0&  1&  0)\\[.1pc]
(1  1  2  3  2  2  1  0) & \myfavouritearrow & ($\phantom{-}$0& $-$1&  0&  1& $-$1&  1&  0& $-$1)\\[.1pc]
(1  2  2  3  2  1  1  0) & \myfavouritearrow & ($\phantom{-}$0&  1&  0&  0&  0& $-$1&  1& $-$1)\\[.1pc]
(1  2  2  3  2  1  1  1) & \myfavouritearrow & ($\phantom{-}$0&  1&  0&  0&  0& $-$1&  0&  1)\\[.1pc]
(1  1  2  3  2  2  1  1) & \myfavouritearrow & ($\phantom{-}$0& $-$1&  0&  1& $-$1&  1& $-$1&  1)\\[.1pc]
(1  1  2  2  2  2  2  1) & \myfavouritearrow & ($\phantom{-}$0&  0&  1& $-$1&  0&  0&  1&  0)\\[.1pc]
(1  1  2  3  3  2  1  0) & \myfavouritearrow & ($\phantom{-}$0& $-$1&  0&  0&  1&  0&  0& $-$1)\\[.1pc]

(1  2  2  3  2  2  1  0) & \myfavouritearrow & ($\phantom{-}$0&  1&  0&  0& $-$1&  1&  0& $-$1)\\[.1pc]
(1  2  2  3  2  2  1  1) & \myfavouritearrow & ($\phantom{-}$0&  1&  0&  0& $-$1&  1& $-$1&  1)\\[.1pc]
(1  1  2  3  3  2  1  1) & \myfavouritearrow & ($\phantom{-}$0& $-$1&  0&  0&  1&  0& $-$1&  1)\\[.1pc]
(1  1  2  3  2  2  2  1) & \myfavouritearrow & ($\phantom{-}$0& $-$1&  0&  1& $-$1&  0&  1&  0)\\[.1pc]
(1  2  2  3  3  2  1  0) & \myfavouritearrow & ($\phantom{-}$0&  1&  0& $-$1&  1&  0&  0& $-$1)\\[.1pc]
(1  2  2  3  3  2  1  1) & \myfavouritearrow & ($\phantom{-}$0&  1&  0& $-$1&  1&  0& $-$1&  1)\\[.1pc]
(1  2  2  3  2  2  2  1) & \myfavouritearrow & ($\phantom{-}$0&  1&  0&  0& $-$1&  0&  1&  0)\\[.1pc]

(1  1  2  3  3  2  2  1) & \myfavouritearrow & ($\phantom{-}$0& $-$1&  0&  0&  1& $-$1&  1&  0)\\[.1pc]
(1  2  2  4  3  2  1  0) & \myfavouritearrow & ($\phantom{-}$0&  0& $-$1&  1&  0&  0&  0& $-$1)\\[.1pc]
(1  2  2  4  3  2  1  1) & \myfavouritearrow & ($\phantom{-}$0&  0& $-$1&  1&  0&  0& $-$1&  1)\\[.1pc]
(1  2  2  3  3  2  2  1) & \myfavouritearrow & ($\phantom{-}$0&  1&  0& $-$1&  1& $-$1&  1&  0)\\[.1pc]
(1  1  2  3  3  3  2  1) & \myfavouritearrow & ($\phantom{-}$0& $-$1&  0&  0&  0&  1&  0&  0)\\[.1pc]
(1  2  3  4  3  2  1  0) & \myfavouritearrow & ($-$1&  0&  1&  0&  0&  0&  0& $-$1)\\[.1pc]
(1  2  3  4  3  2  1  1) & \myfavouritearrow & ($-$1&  0&  1&  0&  0&  0& $-$1&  1)\\[.1pc]
(1  2  2  4  3  2  2  1) & \myfavouritearrow & ($\phantom{-}$0&  0& $-$1&  1&  0& $-$1&  1&  0)\\[.1pc]
(1  2  2  3  3  3  2  1) & \myfavouritearrow & ($\phantom{-}$0&  1&  0& $-$1&  0&  1&  0&  0)\\[.1pc]
(2  2  3  4  3  2  1  0) & \myfavouritearrow & ($\phantom{-}$1&  0&  0&  0&  0&  0&  0& $-$1)\\[.1pc]
(2  2  3  4  3  2  1  1) & \myfavouritearrow & ($\phantom{-}$1&  0&  0&  0&  0&  0& $-$1&  1)\\[.1pc]
(1  2  3  4  3  2  2  1) & \myfavouritearrow & ($-$1&  0&  1&  0&  0& $-$1&  1&  0)\\[.1pc]
(1  2  2  4  3  3  2  1) & \myfavouritearrow & ($\phantom{-}$0&  0& $-$1&  1& $-$1&  1&  0&  0)\\[.1pc]
(2  2  3  4  3  2  2  1) & \myfavouritearrow & ($\phantom{-}$1&  0&  0&  0&  0& $-$1&  1&  0)\\[.1pc]
(1  2  3  4  3  3  2  1) & \myfavouritearrow & ($-$1&  0&  1&  0& $-$1&  1&  0&  0)\\[.1pc]
(1  2  2  4  4  3  2  1) & \myfavouritearrow & ($\phantom{-}$0&  0& $-$1&  0&  1&  0&  0&  0)\\[.1pc]
(2  2  3  4  3  3  2  1) & \myfavouritearrow & ($\phantom{-}$1&  0&  0&  0& $-$1&  1&  0&  0)\\[.1pc]
(1  2  3  4  4  3  2  1) & \myfavouritearrow & ($-$1&  0&  1& $-$1&  1&  0&  0&  0)\\[.1pc]
(1  2  3  5  4  3  2  1) & \myfavouritearrow & ($-$1& $-$1&  0&  1&  0&  0&  0&  0)\\[.1pc]
(2  2  3  4  4  3  2  1) & \myfavouritearrow & ($\phantom{-}$1&  0&  0& $-$1&  1&  0&  0&  0)\\[.1pc]
(1  3  3  5  4  3  2  1) & \myfavouritearrow & ($-$1&  1&  0&  0&  0&  0&  0&  0)\\[.1pc]
(2  2  3  5  4  3  2  1) & \myfavouritearrow & ($\phantom{-}$1& $-$1& $-$1&  1&  0&  0&  0&  0)%\\[.1pc]
\end{tabular}$$

$$\begin{tabular}{lc@{\ \ \ \ \ \ \ \ \ \ }rrrrrrrr}
(2  3  3  5  4  3  2  1) & \myfavouritearrow &($\phantom{-}$1&  1& $-$1&  0&  0&  0&  0&  0)\\[.1pc]
(2  2  4  5  4  3  2  1) & \myfavouritearrow &($\phantom{-}$0& $-$1&  1&  0&  0&  0&  0&  0)\\[.1pc]
(2  3  4  5  4  3  2  1) & \myfavouritearrow &($\phantom{-}$0&  1&  1& $-$1&  0&  0&  0&  0)\\[.1pc]
(2  3  4  6  4  3  2  1) & \myfavouritearrow &($\phantom{-}$0&  0&  0&  1& $-$1&  0&  0&  0)\\
(2  3  4  6  5  3  2  1) & \myfavouritearrow &($\phantom{-}$0&  0&  0&  0&  1& $-$1&  0&  0)\\
(2  3  4  6  5  4  2  1) & \myfavouritearrow &($\phantom{-}$0&  0&  0&  0&  0&  1& $-$1&  0)\\
(2  3  4  6  5  4  3  1) & \myfavouritearrow &($\phantom{-}$0&  0&  0&  0&  0&  0&  1& $-$1)\\
(2  3  4  6  5  4  3  2) & \myfavouritearrow &($\phantom{-}$0&  0&  0&  0&  0&  0&  0&  1)\\
\end{tabular}$$}

\section*{Acknowledgements}

The author would like to thank Prof.~C~S~Seshadri for his constant
encouragement and V~Balaji for the discussions and comments throughout
the preparation of this paper. The referee's comments and suggestions
(in particular the proof of Lemma 8) are gratefully acknowledged. The
author also would like to thank S~Senthamaraikannan and K~V~Subramanian
for carefully reading the preliminary draft of this paper.

\end{document}